\documentclass[onefignum,onetabnum]{siamart190516}

\usepackage{lipsum}
\usepackage{amsfonts}
\usepackage{graphicx}
\usepackage{epstopdf}
\usepackage{algorithmic}
\usepackage{amssymb}
\ifpdf
  \DeclareGraphicsExtensions{.eps,.pdf,.png,.jpg}
\else
  \DeclareGraphicsExtensions{.eps}
\fi

\newsiamremark{remark}{Remark}
\newsiamremark{hypothesis}{Hypothesis}
\crefname{hypothesis}{Hypothesis}{Hypotheses}
\newsiamthm{claim}{Claim}

\headers{Energy stable and MBP-preserving schemes for the $Q$-tensor flow}{D. Hou, X. Li, Z. Qiao and N. Zheng}

\title{Energy stable and maximum bound principle preserving schemes for the $Q$-tensor flow of liquid crystals\thanks{
Submitted to the editors DATE.
\funding{D. Hou is partially supported by Natural Science Foundation of China grant 12001248,  NSF of Jiangsu Province grant BK20201020,  Jiangsu Province Universities  Science Foundation grant 20KJB110013 and  Hong Kong Polytechnic University grant 1-W00D.
X. Li is supported  in part by the National Natural Science Foundation of China  grants  12271302, 12131014.
Z. Qiao is partially supported by Hong Kong Research Council RFS Grant RFS2021-5S03 and GRF Grant 15302122, Hong Kong Polytechnic University Grant 4-ZZLS, and CAS AMSS-PolyU Joint Laboratory of Applied Mathematics.
N. Zheng is partially supported by the Hong Kong Polytechnic University Postdoctoral Research Fund 1-W22P.
}}}

\author{Dianming Hou\thanks{School of Mathematics and Statistics, Jiangsu Normal University, 221116 Xuzhou, China. Current address: Department of Applied Mathematics, The Hong Kong Polytechnic University, Hung Hom, Kowloon, Hong Kong.
  (\email{dmhou@stu.xmu.edu.cn}).}
  \and Xiaoli Li\thanks{School of Mathematics, Shandong University, Jinan, Shandong, 250100, P.R. China.
  (\email{xiaolimath@sdu.edu.cn}).}
\and Zhonghua Qiao\thanks{Department of Applied Mathematics \& Research Institute for Smart Energy, The Hong Kong Polytechnic University, Hung Hom, Kowloon, Hong Kong.
  (\email{zqiao@polyu.edu.hk}).}
  \and Nan Zheng\thanks{Corresponding author. Department of Applied Mathematics, The Hong Kong Polytechnic University, Hung Hom, Kowloon, Hong Kong.
  (\email{nanzheng@polyu.edu.hk}).}
}

\usepackage{amsopn}

\usepackage{subfigure}
\allowdisplaybreaks[4]
\ifpdf
\hypersetup{
  pdftitle={Energy stable and maximum bound principle preserving schemes for the $Q$-tensor flow of liquid crystals},
  pdfauthor={Dianming Hou, Xiaoli Li, Zhonghua Qiao, and Nan Zheng}
}
\fi

\begin{document}
\graphicspath{{figures/},}
\maketitle

\begin{abstract}
In this paper, we propose two efficient fully-discrete schemes for $Q$-tensor flow of liquid crystals by using the first- and second-order stabilized exponential scalar auxiliary variable (sESAV) approach in time and the finite difference method for spatial discretization. The modified discrete energy dissipation laws are unconditionally satisfied for both two constructed schemes. A particular feature is that, for two-dimensional (2D) and a kind of three-dimensional (3D) $Q$-tensor flows, the unconditional maximum-bound-principle (MBP) preservation of the constructed first-order scheme is successfully established, and the proposed second-order scheme preserves the discrete MBP property with a mild restriction on the time-step sizes. Furthermore, we rigorously derive the corresponding error estimates for the fully-discrete second-order schemes by using the built-in stability results.
Finally, various numerical examples validating the theoretical results, such as the orientation of liquid crystal in 2D and 3D, are presented for the constructed schemes.
\end{abstract}


\section{Introduction}
Liquid crystals are generally known as the fourth state of matter \cite{han2015microscopic,lin1995nonparabolic,liu2007dynamics,zhao2016decoupled}, different from gases, liquids, and solids.
They not only  have the property of fluid flows but also possess the property of solids. For example, their molecules have a crystal-like configuration. The orientational order plays an important role in liquid crystals, which is an intermediate state of matter between solid and liquid \cite{ball2011orientability,hu2009stable,huang2016finite,jiang2015stability}.
According to the orientational and the positional order, there are several types of liquid crystal phases, such as the nematic phase, the cholesteric phase and the smectic phase.
The simplest liquid crystal phase is the nematic phase, whose rod-like molecules lack translational order but do exhibit some degree of long-range orientational order.
Within the framework of the Landau–de Gennes theory \cite{de1993physics}, it is widely believed that the local orientation and order of liquid crystal molecules are characterized by a symmetric, traceless $d\times d$ tensor, which is called the $Q$-tensor \cite{han2019transition,sonnet2012dissipative,wang2021modelling}.
We concentrate on the $Q$-tensor theory for nematic liquid crystal flows in this paper.

The Landau-de Gennes free energy functional  \cite{cai2017stable,gudibanda2022convergence,iyer2015dynamic} is given as
\begin{equation*}
	\mathcal{E}[Q]=\int_{\Omega} \mathcal{F}(Q(\bm{x})) d \bm{x},
\end{equation*}
where $\Omega$ is a smooth, bounded domain in $\mathbb{R}^d (d=2, 3)$ and the tensor function $Q$ is defined in the space
\bq\label{eq1}
\mathcal{S}^{(d)} \stackrel{\text { def }}{=}\Big\{M \in \mathbb{R}^{d \times d} \Big| {\color{black}\operatorname{tr}(M)}=0, M^{i j}=M^{j i} \in \mathbb{R},\ \forall i, j=1, \ldots, d\Big\} ,
\eq
{\color{black} where $\operatorname{tr}(M):=\sum_{i=1}^d M^{i i}$.}
Here,
$\mathcal{F}(Q)=\mathcal{F}_{\mathrm{el}}+\mathcal{F}_{\text {bulk }}.$
The elastic part of the free energy density functional $\mathcal{F}$, denoted by $\mathcal{F}_{\mathrm{el}}$, depends on the gradient of $Q$, while the bulk part, denoted by $\mathcal{F}_{\text {bulk }}$, depends exclusively on $Q$.

The bulk-free energy density $\mathcal{F}_{\text {bulk }}$ often takes the form of a truncated expansion. One can use the following simplest setting
\begin{equation}\label{Q4}
	\mathcal{F}_{\text {bulk }} \stackrel{\text { def }}{=} \frac{a}{2} \operatorname{tr}(Q^2)-\frac{b}{3} \operatorname{tr}(Q^3)+\frac{c}{4} \operatorname{tr}^2(Q^2),
\end{equation}
where the constants for the bulk material are taken to be $a \in \mathbb{R}$,\ $b\geq0$ and $c>0$. The phase, which changes from the isotropic phase to the nematic phase, is controlled by $a$ \cite{zhao2017novel}. Specifically,  the system is in the
isotropic phase with $a>0$ while it is in the nematic phase with $a<0$.  It can be obtained that $b\geq0$ in \cite{iyer2015dynamic,wu2019dynamics}, here  $c>0$ \cite{majumdar2010equilibrium,wang2023q} ensures that the $\mathcal{F}_{\text {bulk}}$  is bounded from below.

The strain energy density is provided by the elastic-free energy density $\mathcal{F}_{\text {el}}$ as a result of spatial variations in the tensor order parameter.
The simplest form which is invariant under rigid rotation and material symmetry is as follows
\begin{equation*}
\mathcal{F}_{\text {el}} \stackrel{\text { def }}{=} \frac{L_1}{2}|\nabla Q|^2+
\sum_{i,j,k=1}^{d}\Big(\frac{L_2}{2} \partial_j Q^{i k} \partial_k Q^{i j}+\frac{L_3}{2} \partial_j Q^{i j} \partial_k Q^{i k}\Big),
\end{equation*}
where $\partial_k Q^{i j}$ represents a partial derivative of $i j$ component of $Q$ with respect to $x_k$.


The following $L^2$-gradient flow in $\mathbb{R}^d$, $d=2,3$, corresponds to the energy functional $\mathcal{E}[Q]$ where $Q$ takes values in $\mathcal{S}^{(d)}$ :
\begin{equation}\label{Q6}
\frac{\partial Q^{i j}}{\partial t}=-\Big(\frac{\delta \mathcal{E}}{\delta Q}\Big)^{i j}+\lambda \delta^{i j}+\mu^{i  j}-\mu^{j i}, \quad 1 \leq i, j \leq d,
\end{equation}
where $\lambda$ {\color{black} and $\mu=(\mu^{i  j})_{d\times d}$  are} Lagrange multipliers corresponding to the tracelessness constraint {\color{black} and  the matrix symmetry constraint, respectively.} The variational derivative of  $\mathcal{E}$  with respect to $Q$ is denoted as $\delta \mathcal{E}/{\delta Q}$.
We always impose $L_i>0$, which ensures that the summation of the first three quadratic terms concerning $L_1, L_2$ and $L_3$ in $\mathcal{F}_{\mathrm{el}}$ is positive definite from a modeling perspective \cite{wang2023q}.


The evolution equation \eqref{Q6} {\color{black} in nonequilibrium situation \cite{gudibanda2022convergence}} is given below after expansion. {\color{black}For more details, please refer to Appendix A in \cite{cai2017stable,iyer2015dynamic}.}
\begin{equation}\label{Q9}
\begin{aligned}
\partial_t Q^{i j}=&\dps L_1 \Delta Q^{i j}
+\frac{L_2+L_3}{2}\Big( \sum_{k=1}^{d}\big(\partial_j \partial_k Q^{i k}
+\partial_i \partial_k Q^{j k}\big)
-\frac{2}{d} \sum_{k,l=1}^{d}\partial_l \partial_k Q^{l k}\delta^{i j}\Big) \\
&\dps-\Big(a Q^{i j}-b\big((Q^2)^{i j}-\frac{1}{d}tr(Q^2)\delta^{ij}\big)+c \operatorname{tr}(Q^2) Q^{i j}\Big), \quad1 \leq i, j \leq d.
\end{aligned}
\end{equation}
   {\color{black}
 	The system is subject to either Dirichlet or periodic boundary conditions, with an initial condition $Q(\bm{x}, 0)=Q^0(\bm{x})$. 
 	The primary objective of this study is to analyze the behavior of nonequilibrium systems.
}

There exist several numerical researches for the $Q$-tensor model \cite{badia2011finite,cruz2013numerical,cai2017stable,gudibanda2022convergence,mori1999multidimensional,tovkach2017q,yang2013modeling,zhao2017novel}.
Cai, Shen and Xu \cite{cai2017stable} proposed a first-order nonlinear scheme for a 2D dynamic $Q$-tensor model of nematic liquid crystals by using a stabilizing technique, the maximum principle and convergence analysis are established. The constructed scheme in this paper can not guarantee unconditional energy stability.
Gudibanda  et al. \cite{gudibanda2022convergence} presented a fully discrete convergent finite difference scheme for the 3D $Q$-tensor flow of liquid crystals based on the IEQ method. They proved the stability properties of the scheme and showed that it converges to weak solutions under a natural assumption. However, it seems that  the maximum bound principle (MBP) can not be satisfied for their proposed schemes. 
Therefore, there is a high demand to construct efficient numerical schemes that can satisfy the unconditional energy dissipation law while preferably preserving the MBP with rigorous error estimates.

The main purpose of this paper is to construct efficient fully-discrete schemes for the $Q$-tensor flow using the stabilized exponential scalar auxiliary variable (sESAV) approach \cite{ju2022generalized,ju2022stabilized,shen2019new}.
The main purposes of this paper are
\begin{itemize}
\item  to construct two linear efficient fully discrete  first- and second-order schemes for the $Q$-tensor flow of liquid crystals.

\item  to establish unconditional energy stability and prove the MBP for both two constructed schemes. 
with an appropriate stabilizing parameter, the first-order scheme preserves the MBP unconditionally, while the second-order scheme does so with mild time-step size constraints.

\item  to carry out a rigorous error analysis for the fully discrete schemes, which is derived by establishing uniform bound for the numerical solution.
\end{itemize}
To the best of the authors' knowledge, this is the first linear, energy-stable and MBP-preserving numerical scheme for the $Q$-tensor flow in both 2D and 3D cases.
Moreover, we also establish the rigorous error analysis for the proposed schemes by using the built-in stability results and a sequence of  auxiliary estimates.

The remainder of the paper is structured as follows.
Section 2 is devoted to the physical properties of the model and some preliminaries.
In Section 3, we construct the fully-discrete first- and second-order scheme respectively, together with the energy dissipation law and MBP.  
In Section 4, we establish the rigorous error estimate for the numerical scheme.
In Section 5, some numerical experiments are carried out to demonstrate the performance of the proposed schemes.
Concluding remarks are given in Section 6.


\section{The physical properties of the model and some preliminaries}
\textcolor{black}{Without loss of generality, we only consider the problem \eqref{Q9} with the  homogenous Dirichlet boundary condition, i.e., $Q|_{\partial \Omega}=0$. It is easy to extend the corresponding obtained results to the cases with the inhomogeneous Dirichlet or periodic boundary conditions. }

Firstly, we provide some notations. For any $\alpha\in\mathbb{R}$, $\alpha_{+}$ and $\alpha_{-}$ are given by
\bq\label{equa2}
\alpha_{+}:=\max\{0,\alpha\},\quad \alpha_{-}:=\max\{0,-\alpha\}.
\eq
For any tensor functions $\Phi, \Psi \in \mathbb{R}^{d\times d}$, we define the Frobenius product as follows:
\begin{equation*}
\Phi: \Psi \stackrel{\text { def }}{=} \operatorname{tr}(\Phi^T \Psi)=\sum_{i,j=1}^{d}\Phi^{i,j}\Psi^{i,j}.
\end{equation*}
The tensor function $\Phi \in \mathbb{R}^{d \times d}$, we use $|\Phi|_{F}$ to denote its Frobenius norm, i.e.,
\begin{equation*}
|\Phi|_F=\sqrt{\operatorname{tr}(\Phi^T \Phi)}.
\end{equation*}
Besides, we define $L^p(1 \leq p \leq \infty)$ tensor-function space by
\begin{equation*}
L^p(\Omega \rightarrow \mathbb{R}^{d \times d}) \stackrel{\text { def }}{=}\big\{\Phi: \Omega \rightarrow \mathbb{R}^{d \times d},|\Phi|_{F} \in L^p(\Omega)\big\}.
\end{equation*}
The corresponding $L^{\infty}$-norm for $\Phi\in L^p(\Omega \rightarrow \mathbb{R}^{d \times d}) $ is given by
\bq\label{eq6}
\|\Phi\|_{L^{\infty}(\Omega)}=\sup_{\bm{x}\in \Omega}|\Phi(\bm{x})|_{F}.
\eq

The gradient flow \eqref{Q6} satisfies the following energy dissipation law:
\begin{equation*}
\frac{d}{d t} \mathcal{E}[Q]=-\int_{\Omega}\Big|\frac{\delta \mathcal{E}}{\delta Q}-\lambda \mathbf{I}_d+\mu-\mu^T\Big|^2 d \bm{x}.
\end{equation*}
Moreover, it has the nature of MBP property in the sense of the Frobenius norm, as stated in the following lemma,
\textcolor{black}{which was given by  Proposition 2.2 in \cite{iyer2015dynamic} and Theorem 1 in \cite{majumdar2010equilibrium}.
	}
\begin{lem}\label{lemma1}
Consider the evolution problem \eqref{Q6} with $a,\; b\in\mathbb{R}$ and $c>0$ on a bounded, smooth domain $\Omega \subset \mathbb{R}^d$. For smooth solutions $Q$, there exists a positive number $\eta$ such that, when $d=2$ or $L_2+L_3=0$ for $d=3$, if $\|Q^0\|_{L^{\infty}(\Omega)} \leq \eta$, then it holds for any $T > 0$ and $t \in [0,T]$ that $\|Q(t)\|_{L^{\infty}(\Omega)} \leq \eta.$
\end{lem}

\begin{proof}
\textcolor{black}{We first establish a priori $L^\infty$-bound for the solution $Q$.}
Due to the fact that $c > 0$, there exists $\eta>0$ such that for each $Q\in S^{(d)}$,
\begin{equation*}
-a\operatorname{tr}(Q^2)+b \operatorname{tr}(Q^3)-c\operatorname{tr}^2(Q^2) \leq 0, \quad \text{if} \
|Q|_{\textcolor{black}{F}}\geq \eta.
\end{equation*}
We assume  $\|Q_0\|_{L^{\infty}(\Omega)}\leq \eta$, multiplying
\eqref{Q9} by $Q(|Q|_{F}^2-\eta^2)_+$ and integrating over $\Omega$,
after integration by parts, we obtain
\begin{equation*}
\begin{aligned}
&\frac{1}{4}\frac{d}{dt}\int_\Omega(|Q|_{F}^2-\eta^2)_+^2d\bm{x}\\
=&-L\int_\Omega |\nabla Q|^2 (|Q|_{F}^2-\eta^2)_+d\bm{x}
-\frac{L}{2}\int_\Omega|\nabla (|Q|_{F}^2-\eta^2)_+|^2d\bm{x}\\
&+\int_\Omega (-a\operatorname{tr}(Q^2)+b \operatorname{tr}(Q^3)-c\operatorname{tr}^2(Q^2))(|Q|_{\textcolor{black}{F}}^2-\eta^2)_+d\bm{x}
\leq 0,
\end{aligned}
\end{equation*}
which implies
$\|Q(\cdot,t)\|_{L^{\infty}(\Omega)}\leq \eta,\quad \forall t\in [0,T]$.
\end{proof}

\begin{remark}\label{remk1}
For the 2D case, $Q\in \mathcal{S}^{(2)} $ implies $\operatorname{tr}(Q^3)=0$, then it follows from \eqref{Q4} and \eqref{Q9} that $b=0$. Furthermore, we obtain the minimum value of $\eta$ in Lemma \ref{lemma1} for $d=2$, as follows:
\bq\label{equ1-1}
\tps\eta^{(2)}=\max\Big\{\|Q^{0}\|_{L^{\infty}(\Omega)},\; \sqrt{\frac{a_{-}}{c}}\Big\},
\eq
where $a_{-}$ is defined in \eqref{equa2}.
As for the 3D case, we recall that $|\operatorname{tr}(Q^3)|\leq |Q|^3_{F}/\sqrt{6}$, which implies the minimum value of $\eta$ for $d=3$: (see e.g. \cite{wang2017q,contreras2019elementary})
\bq\label{equ1}
\eta^{(3)}=
\begin{cases}
\begin{array}{r@{}l}
&\max\Big\{\|Q^{0}\|_{L^{\infty}(\Omega)},\; \frac{|b|+\sqrt{b^2-24ac}}{2\sqrt{6}c}\Big\},  \qquad\mbox{ if }  a\leq\frac{b^{2}}{24c},\\
&\|Q^{0}\|_{L^{\infty}(\Omega)}, \qquad \qquad \qquad  \qquad \qquad\qquad \mbox{ otherwise.}
\end{array}
\end{cases}
\eq
From the definitions of \eqref{equ1-1} and \eqref{equ1}, it also implies that $\|Q^{0}\|_{L^{\infty}(\Omega)}\leq \eta^{(d)},\; d=2,3.$
\end{remark}

\section{Fully discrete sESAV schemes and structure-preserving properties}
In this section, we construct the first- and second-order structure-preserving fully discrete schemes for the $Q$-tensor flow of liquid crystals  using finite difference method.

Recalling Lemma \ref{lemma1} and noting the fact that $\mathcal{F}_{\text {bulk}} $ is bounded from below \textcolor{black}{\cite{zhao2017novel}}.
\textcolor{black}{
There exist two  positive constants $C_\ast, C^\ast$
such that}
\bq\label{equa3}
-C_\ast\leq\mathcal{E}_1[Q(t)]=\int_{\Omega}\frac{a}{2} \operatorname{tr}\left(Q^2\right)
-\frac{b}{3} \operatorname{tr}\left(Q^3\right)
+ \frac{c}{4} \operatorname{tr}^2\left(Q^2\right)
d \bm{x}\leq C^\ast,
\eq
{\color{black} where $C_\ast$ only depends on $a,b, c$, and can be determined using Theorem 2.1 in \cite{zhao2017novel}.}

Introducing $s(t)=\mathcal{E}_1[Q(t)]$ as an exponential scalar auxiliary variable,
we then have the following energy form
\brr\label{rev_eq1}
E[Q,s]&\dps=\int_{\Omega}
\Big[\frac{L_1}{2}|\nabla Q|^2+\frac{L_2+L_3}{2}\sum_{i,j,k=1}^{d} \partial_j Q^{i j} \partial_k Q^{i k}\Big]d\bm{x}
+s\\
&\textcolor{black}{\dps=\frac{L_1}{2}\|\nabla Q\|^2+
\frac{L_2+L_3}{2}\|\text{div} Q\|^{2}+s}\\
&{\color{black}\geq -C_\ast}.
\err
\textcolor{black}{For the nonequilibrium system considered in this paper,
it has been stated in \cite{gudibanda2022convergence} that $E[Q,s]$ is equivalent to  $\mathcal{E}[Q]$.}
We rewrite \eqref{Q9} as the following \cite{ju2022stabilized}:
\begin{subequations}\label{Q20}
\begin{align}
\partial_t Q^{i j}
=&L_1 \Delta Q^{i j}
+\frac{L_2+L_3}{2}\Big(\sum_{k=1}^{d} \big(\partial_j \partial_k Q^{i k}
+\partial_i \partial_k Q^{j k}\big)
- \frac{2}{d}\sum_{k,l=1}^{d}\partial_l \partial_k Q^{l k}\delta^{i j}\Big) \\
&+\frac{\exp\{s\}}{\exp\{\mathcal{E}_1[Q]\}}f(Q^{ij}),\nonumber\\
s_t=&-\frac{\exp\{s\}}{\exp\{\mathcal{E}_1[Q]\}}\sum_{i,j=1}^{d}	\big(f(Q^{ij}), Q^{i j}_t
\big),
\end{align}
\end{subequations}
where $f(Q^{ij}):=-a Q^{i j}+b\big((Q^2)^{i j}-\frac{1}{d}tr(Q^2)\delta^{ij}\big)-c \operatorname{tr}\left(Q^2\right) Q^{i j}.$


%

\subsection{Fully discrete sESAV finite difference schemes}
To simplify the notation, we consider
the problem \eqref{Q20} in the domain $\Omega = [0, L_d]^3$.
Given a positive integer $M$, the uniform mesh partitioning size for each spatial direction is set to be $h=L_d/M$. We denote by $\mathbf{E}$ and $\mathbf{C}$ the two sets of mesh points, defined by
$$
\mathbf{E}=\big\{x_{p}=ph\;\big|\; p=0,1, \cdots,M\big\},\qquad \mathbf{C}=\big\{x_{p+\frac{1}{2}}=\big(p+\textstyle\frac{1}{2}\big)h\;\big|\; p=0,1, \cdots,M-1\big\}.
$$
The necessary grid function spaces are given by
\bry
\mathbf{E}^{0}_{h}=&\dps\big\{U: \mathbf{E}\times\mathbf{E}\times\mathbf{E}\rightarrow\mathbb{R}\;\big|\;U_{p,q,r},\; U_{0,q,r}=U_{p,0,r}=U_{p,q,0}=0, \;0\leq p,q,r\leq M\big\},\\[4pt]
e^{(1)}_{h}=&\dps\big\{U: \mathbf{C}\times\mathbf{E}\times\mathbf{E}\rightarrow\mathbb{R}\;\big|\;U_{p+\frac{1}{2},q,r}, \;0\leq p\leq M-1,\;0\leq q,r\leq M\big\},\\[4pt]
e^{(2)}_{h}=&\dps\big\{U: \mathbf{E}\times\mathbf{C}\times\mathbf{E}\rightarrow\mathbb{R}\;\big|\;U_{p,q+\frac{1}{2},r}, \;0\leq q\leq M-1,\;0\leq p,r\leq M\big\},\\[4pt]
e^{(3)}_{h}=&\dps\big\{U: \mathbf{E}\times\mathbf{E}\times\mathbf{C}\rightarrow\mathbb{R}\;\big|\;U_{p,q,r+\frac{1}{2}}, \;{\color{black}0\leq r\leq M-1},\; 0\leq p,q\leq M\big\}.
\ery
For ease of notation, we will also impose homogenous Dirichlet boundary condition on the ghost nodes, that is, for any $-1\leq p,q,r\leq M+1,$
$$U_{-1,q,r}=U_{M+1,q,r}=U_{p,-1,r}=U_{p,M+1,r}=U_{p,q,-1}=U_{p,q,M+1}=0.$$
Then, we give several difference operators for the grid function $\{U_{p,q,r}\}_{p,q,r=0}^{M}$:
\bryl
\big(D_{1}^{+}U\big)_{p+\frac{1}{2},q,r}=\frac{U_{p+1,q,r}-U_{p,q,r}}{h},~~\big(D_{2}^{+}U\big)_{p,q+\frac{1}{2},r}=\frac{U_{p,q+1,r}-U_{p,q,r}}{h}, \\
\big(D_{3}^{+}U\big)_{p,q,r+\frac{1}{2}}=\frac{U_{p,q,r+1}-U_{p,q,r}}{h},~~\big(D_{1}^{-}U\big)_{p,q,r}=\frac{U_{p+\frac{1}{2},q,r}-U_{p-\frac{1}{2},q,r}}{h},\\
\big(D_{2}^{-}U\big)_{p,q,r}=\frac{U_{p,q+\frac{1}{2},r}-U_{p,q-\frac{1}{2},r}}{h},~~ \big(D_{3}^{-}U\big)_{p,q,r}=\frac{U_{p,q,r+\frac{1}{2}}-U_{p,q,r-\frac{1}{2}}}{h}\\
\big(D_{1}^{c}U\big)_{p,q,r}=\frac{U_{p+1,q,r}-U_{p-1,q,r}}{2h}, ~~\big(D_{2}^{c}U\big)_{p,q,r}=\frac{U_{p,q+1,r}-U_{p,q-1,r}}{2h},\\
\big(D_{3}^{c}U\big)_{p,q,r}=\frac{U_{p,q,r+1}-U_{p,q,r-1}}{2h}.\\
\eryl
Moreover, we denote by $D^{c}_{k,l}U=D^{c}_{k}D^{c}_{l}U,\; k,l=1,2,3$ for simplicity.
The discrete gradient operator $\nabla_{h}$: $E^{0}_{h}\rightarrow(e^{(1)}_{h},e^{(2)}_{h},e^{(3)}_{h})^{T}$ is given by
\beq
\big(\nabla_{h}U\big)_{p,q,r}=\Big(\big(D^{+}_{1}U\big)_{p+\frac{1}{2},q,r},\big(D^{+}_{2}U\big)_{p,q+\frac{1}{2},r},\big(D^{+}_{3}U\big)_{p,q,r+\frac{1}{2}}\Big)^{T},
\eeq
and the discrete divergence operator $\nabla_{h}\cdot$: $(e^{(1)}_{h},e^{(2)}_{h},e^{(3)}_{h})^{T}\rightarrow E^{0}_{h}$ is represented by
\beq
\nabla_{h}\cdot (U^{(1)},U^{(2)},U^{(3)})^{T}=D^{-}_{1}U^{(1)}+D^{-}_{2}U^{(2)}+D^{-}_{3}U^{(3)},
\eeq
where $(U^{(1)},U^{(2)},U^{(3)})^{T}\in (e^{(1)}_{h},e^{(2)}_{h},e^{(3)}_{h}).$ Therefore, the discrete Laplacian $\Delta_{h}$: $E^{0}_{h}\rightarrow E^{0}_{h}$ is defined by $\Delta_{h}U=\nabla_{h}\cdot(\nabla_{h}U)$
for any grid function $U\in E^{0}_{h}.$
We also introduce several discrete inner products:
\bry
&\big<U,V\big>_{h}=\dps h^{3}\sum_{p,q,r=1}^{M-1}U_{p,q,r}V_{p,q,r},\quad \forall\,U,V\in E^{0}_{h},\\[7pt]
&[U^{(k)},V^{(k)}]_{k,h}=\big<a_{k}(U^{(k)}V^{(k)}),1\big>_{h},\quad \forall\,U^{(k)},V^{(k)}\in e^{(k)}_{h},\ \ k=1,2,3,\\[7pt]
&\dps [(U^{(1)},U^{(2)},U^{(3)})^{T},(V^{(1)},V^{(2)},V^{(3)})^{T}]_{h}=\sum_{k=1}^{3}[U^{(k)},V^{(k)}]_{k,h},
\ery
where $a_{k}: e^{k}_{h}\rightarrow E^{0}_h, k=1,2,3$ are three average operators, defined by
\bryl
(a_{1}U^{(1)})_{p,q,r}=\dps\frac12({U^{(1)}_{p+\frac{1}{2},q,r}+U^{(1)}_{p-\frac{1}{2},q,r}}),\ \ (a_{2}U^{(2)})_{p,q,r}=\dps\frac12({U^{(2)}_{p,q+\frac{1}{2},r}+U^{(2)}_{p,q-\frac{1}{2}}},r),\\
(a_{3}U^{(3)})_{p,q,r}=\dps\frac12({U^{(3)}_{p,q,r+\frac{1}{2}}+U^{(3)}_{p,q,r-\frac{1}{2}}}),\quad 0\leq p,q,r\leq M,
\eryl
for any $U^{(k)}\in e^{k}_{h},\: k=1,2,3.$ Then, for any $U\in \mathbf{E}^{0}_{h},$ the corresponding discrete $L^{2}$, $H^{1}$ and $L^{\infty}$ norms are given by
\bry
&\dps\|U\|^{2}_{h}=\big<U,U\big>_{h}, \ \ \|\nabla_{h}U\|^{2}_{h}=[\nabla_{h}U,\nabla_{h}U]_{h}=\sum_{k=1}^{3}[D^{+}_{k}U,D^{+}_{k}U]_{k,h},\\
&\dps\|U\|^{2}_{H^{1}_{h}}=\|U\|^{2}_{h}+\|\nabla_{h}U\|^{2}_{h},\quad \dps{\color{black}\|\vec{U}\|_{\infty}=\max_{0\leq p,q,r\leq M}|U_{p,q,r}|,}
\ery
{\color{black}where $\vec{U}$ represents $U$ 
as a vector with its components arranged along the directions $x_1,\cdots,x_d$.
	}
Using the summation-by-parts, it is easy to check that, for any $U,V\in \mathbf{E}^{0}_{h},$ it holds
\brl\label{eq2}
\big<D^{-}_{k}D^{+}_{k}U,V\big>_{h}\dps=-\big<D^{+}_{k}U,D^{+}_{k}V\big>_{h}=\big<U,D^{-}_{k}D^{+}_{l}V\big>_{h}, \\
\big<D^{c}_{k,l}U,V\big>_{h}\dps=-\big<D^{c}_{l}U,D^{c}_{k}V\big>_{h}=-\big<D^{c}_{k}U,D^{c}_{l}V\big>_{h}=\big<U,D^{c}_{k,l}V\big>_{h}, \\
\erl
which also imposes that 
\beq
\big<\Delta_h U,V\big>_{h}=-\big[\nabla_h U, \nabla_h V\big]=\big<U,\Delta_h V\big>_{h}.
\eeq

Analogous to the tensor-function space $\mathcal{S}^{(d)}$ in \eqref{eq1}, we also define grid tensor-function space satisfying the homogenous Dirichlet boundary condition for $d=3$,
$$\mathbf{S}^{(3)}_{h}=\Big\{\Phi_{h}\Big| \Phi^{i,j}_{h}\in\mathbf{E}^{0}_{h},\; \mbox{tr}(\Phi_{h}):=\sum_{i=1}^{3}\Phi^{ii}_{h}=0,\; \Phi^{i,j}_{h}=\Phi^{j,i}_{h},\; i,j=1,2,3\Big\}.$$
We can give a similar definition of general tensor-function space $\mathbf{S}^{d}$ for $d=2.$ Then, the corresponding discrete $L^{2}$ norm and the discrete $H^{1}$ semi-norm of the grid tensor-functions $\Phi_{h},\Psi_{h}\in \mathbf{S}^{(3)}_{h}$, are defined respectively by
\beq
 \dps \|\Phi_{h}\|_{h}:=\sqrt{\sum_{i,j=1}^{3}\|\Phi^{i,j}_{h}\|^{2}_{h}},\;\quad
\|\nabla_{h}{\color{black}\Psi_{h}}\|_{h}:=\sqrt{\sum_{i,j=1}^{3}\|\nabla_{h}{\color{black}\Psi^{i,j}_{h}}\|^{2}_{h}}.
\eeq
Moreover, we define following tensor-matrix-product and discrete Frobenius-product for any $\Phi_{h},\Psi_{h}\in \mathbf{S}^{(3)}_{h}$, respectively, by
\beq
\dps \big(\Phi_{h}\Psi_{h}\big)^{i,j}:=\sum_{k=1}^{3}\Phi^{i,k}_{h}.*\Psi^{k,j}_{h}, \; i,j=1,2,3; \ \ \Phi_{h}:\Psi_{h}:=\mbox{tr}(\Phi_{h}\Psi_{h})=\sum_{i,j=1}^{3}\Phi^{i,j}_{h}.*\Psi^{i,j}_{h},\\
\eeq
where $\Phi^{i,k}_{h}.*\Psi^{k,j}_{h}$ represents the point-wise multiplication of matrices, that is $$\big(\Phi^{i,k}_{h}.*\Psi^{k,j}_{h}\big)_{p,q,r}=\big(\Phi^{i,k}_{h}\big)_{p,q,r}\cdot\big(\Psi^{k,j}_{h}\big)_{p,q,r},\ \ \ p,q,r=0,1,\cdots,M.$$

Now, we are in a position to construct fully discrete first- and second-order in-time schemes for the Q-tensor model \eqref{Q6}.
The first-order stabilized ESAV (sESAV1) fully-discrete scheme  for \eqref{Q20} reads:
given the initial values of $(Q_{h}^0, s^0)$, then for
any $n \geq 0$, we find $(Q_{h}^{n+1}, s^{n+1})\in\mathbf{S}^{(3)}_{h}\times\mathbb{R}$ such that for  $i,j=1,2,3$
\begin{subequations}\label{Q23}	
\begin{align}
\frac{Q_{h}^{i j, n+1}-Q_{h}^{i j, n}}{\tau}
=& L_1 \Delta_h Q_{h}^{i j, n+1}+\frac{L_2+L_3}{2}\mathcal{D}^{c}_{h}Q_{h}^{i j, n+1}\nonumber\\
&+g(Q_{h}^n, s^n)\big(f(Q_{h}^{ij,n})-\kappa(Q_{h}^{i j, n+1}-Q_{h}^{i j, n})\big), \label{Q23a}\\
\frac{\textcolor{black}{\tilde{s}}^{n+1}-s^{n}}{\tau}=&-g(Q_{h}^n, s^n)\sum_{i,j=1}^{3}\Big<f(Q_{h}^{ij, n}),	\frac{Q_{h}^{i j, n+1}-Q_{h}^{i j, n}}{\tau}\Big>_{h},\label{Q23b}\\
\textcolor{black}{s^{n+1}=}&\textcolor{black}{\max\left\{\tilde{s}^{n+1},-C_\ast-E_{el}[Q^{n+1}_{h}]\right\}}, \label{Q23c}
\end{align}
\end{subequations}
where $\kappa\geq 0$ is a stabilizing constant, $g(Q^{n}_h, s^{n}_h):=\exp\{s_h\}/\exp\{\mathcal{E}_{1h}[Q_h]\} >0$, and
$$
\textcolor{black}{E_{el}[Q^{n+1}_{h}]=\frac{L_1}{2}\|\nabla_h Q^{n+1}_h\|_h^2+
\frac{L_2+L_3}{2}\sum_{i=1}^{3} \|D^{c}_1 Q^{i 1, n+1}_{h}+ D^{c}_2 Q^{i 2, n+1}_{h}+D^{c}_3 Q^{i 3, n+1}_{h}\|^{2}_{h}.}
$$
Here, the operator $\mathcal{D}^{c}_{h} Q_{h}^{i j, n+1}$ is given by
\bq\label{optor}
\mathcal{D}^{c}_{h} Q_{h}^{i j, n+1}:=\sum_{k=1}^{3}\big(D_{jk}^{c} Q_{h}^{i k, n+1}
+D_{ik}^{c} Q_{h}^{j k, n+1}\big)-\frac{2}{3} \sum_{k,l=1}^{3}D_{lk}^{c} Q_{h}^{l k, n+1}\delta^{i j}.
\eq
We denote the above first-order scheme as $[Q^{n+1}_{h},s^{n+1}]=\mbox{sESAV1}(Q^{n}_{h},s^{n},\tau).$
We point out that the constructed scheme  \eqref{Q23} preserves the trace-free and symmetry property of $Q$, which can be obtained by using a similar proof in \cite{gudibanda2022convergence}.
\begin{pro}\label{pro1}
If $Q^n_{h}$ is trace-free and symmetric, then $Q_{h}^{n+1}$ computed
by the sESAV1 scheme \eqref{Q23} is also trace-free and symmetric.
\end{pro}

The unconditionally energy stability of the sESAV1 scheme \eqref{Q23} is stated in the following theorem.
\begin{thm} \label{thm2}
For any $\kappa\geq 0$, the \textnormal{sESAV1} scheme \eqref{Q23} is unconditionally energy stable in the sense that
$$E_h[Q_{h}^{n+1},s^{n+1}]\leq E_h[Q_{h}^n,s^n],\quad n\geq0.$$
\end{thm}

\begin{proof}
{\color{black}We apply the mathematical induction method to complete this proof.  Since $s^{0}+E_{el}[Q^{0}_{h}]=\mathcal{E}_1[Q^{0}_{h}]+E_{el}[Q^{0}_{h}]\geq -C_*$,
 we assume $s^{k}+E_{el}[Q^{k}_{h}]\geq -C_*$ for all $0\leq k\leq n$.}
Taking the discrete inner product with \eqref{Q23a} by $Q_{h}^{ij,n+1}-Q_h^{ij,n}$, using the discrete integral by parts \eqref{eq2} and summing up $i,j$ from 1 to 3, yields
\begin{equation}\label{Q37}
\begin{aligned}
&(\frac{1}{\tau}
+\kappa g(Q_{h}^n, s^n))
\|Q_{h}^{n+1}-Q_{h}^n\|^2_{h}\\
\leq&-\frac{L_1}{2}\big(\|\nabla_hQ_{h}^{n+1}\|_{h}^2-\|\nabla_h Q_h^{n}\|_{h}^2\big)+\textcolor{black}{\frac{L_2+L_3}{2}}\sum_{i,j=1}^{3}\big<\mathcal{D}^{c}_{h} Q_{h}^{i j, n+1}, Q_{\textcolor{black}{h}}^{i \textcolor{black}{j}, n+1}\\
&-Q_{\textcolor{black}{h}}^{i \textcolor{black}{j}, n}\big>_{h}+g(Q_{h}^n, s^n)
\sum_{i,j=1}^{3}\big<f(Q^{ij,n}), Q_{\textcolor{black}{h}}^{i j, n+1}-Q_{\textcolor{black}{h}}^{i j, n}\big>_{h}.
\end{aligned}
\end{equation}
Moreover, we use \eqref{eq2} and $Q^{ij}_{h}=Q^{ji}_{h}$ for $i,j=1,2,3$, to obtain that
\brr	\label{Q38}
&\dps\sum_{i,j=1}^{3}\Big<\sum_{k=1}^{3}\big(D_{jk}^{c} Q_{h}^{i k, n+1}
+D_{ik}^{c} Q_{h}^{j k, n+1}\big), Q_{\textcolor{black}{h}}^{i j, n+1}-Q_{\textcolor{black}{h}}^{i j, n}\Big>_{h}\\
=&\dps-2\sum_{i=1}^{3}\Big<\sum_{k=1}^{3}D_{k}^{c} Q_{h}^{i k, n+1},\sum_{k=1}^{3}D_{k}^{c}\big(Q_{\textcolor{black}{h}}^{i k, n+1}-Q_{\textcolor{black}{h}}^{i k, n}\big)\Big>_{h}\\
\leq&\dps-\sum_{i=1}^{3}\Big[\Big\|\sum_{k=1}^{3}D_{k}^{c} Q_{h}^{i k, n+1}\Big\|^{2}_{h}-\Big\|\sum_{k=1}^{3}D_{k}^{c} Q_{h}^{i k, n}\Big\|^{2}_{h}\Big].
\err
In addition, using the trace-free preservation of the sESAV1 scheme \eqref{Q23}, we derive
\bq\label{eq3}
\sum_{i,k,l=1}^{3}\big<D_{lk}^{c} Q_{h}^{l k, n+1}, Q_{h}^{i i, n+1}-Q_{h}^{i i, n}\big>_{h}=0
\eq	
It follows from  \eqref{Q23b}-\eqref{Q38} and \eqref{eq3} that
\beq
E_h[Q_{h}^{n+1},\textcolor{black}{\tilde{s}}^{n+1}]-E_h[Q_h^{n},s^{n}]=\dps-\Big(\frac{1}{\tau}
+\kappa g(Q_{h}^n, s^n)\Big)\|Q_{h}^{n+1}-Q_{h}^n\|_{h}^2
\leq0.
\eeq
\textcolor{black}{If $\tilde{s}^{n+1}\geq-C_\ast-E_{el}[Q^{n+1}_{h}]$, it follows that $s^{n+1}=\tilde{s}^{n+1}$, and, consequently,
\begin{equation}\label{E_h}
	E_h[Q_{h}^{n+1},s^{n+1}]-E_h[Q_h^{n},s^{n}]\leq 0.
\end{equation}
In the case where $\tilde{s}^{n+1} < -C_\ast-E_{el}[Q^{n+1}_{h}]$, it can be observed that
$$-C_\ast=E_h[Q_{h}^{n+1},s^{n+1}]
\leq s^{n}+E_{el}[Q^{n}_{h}]=E_h[Q_h^{n},s^{n}],$$
resulting in the direct derivation of  \eqref{E_h}.
}
\end{proof}

The second-order in time stabilized ESAV (sESAV2) fully-discrete scheme  for the problem \eqref{Q20} is presented as follows: given the initial values of ($Q_{h}^0$, $s^0$), then for
any $n \geq 0$, we obtain $(Q_{h}^{n+1}, s^{n+1})\in \mathbf{S}^{3}_{h}\times\mathbb{R}$ through the following
\begin{subequations}\label{Q43}
\begin{align}
\frac{Q_{h}^{i j, n+1}-Q_{h}^{i j, n}}{\tau}
=&L_1 \Delta_h Q_{h}^{i j, n+\frac{1}{2}}
+\textcolor{black}{\frac{L_2+L_3}{2}}\mathcal{D}^{c}_{h}Q^{ij,n+\frac{1}{2}}_{h}\nonumber\\
&+g(Q_{*,h}^{n+\frac{1}{2}}, s_{*}^{n+\frac{1}{2}})
\big(f(Q_{*,h}^{ij,n+\frac{1}{2}})
-\kappa (Q_{h}^{i j, n+\frac{1}{2}}-Q_{*,h}^{ij,n+\frac{1}{2}})\big),\label{Q43a}\\
\frac{\textcolor{black}{\tilde{s}}^{n+1}-s^{n}}{\tau}=&-	g(Q_{*,h}^{n+\frac{1}{2}}, s_{*}^{n+\frac{1}{2}})\sum_{i,j=1}^{3}\Big<f(Q_{*,h}^{ij,n+\frac{1}{2}}), 	 \frac{Q_{h}^{i j, n+1}-Q_{h}^{i j, n}}{\tau}
\Big>_{h}\nonumber\\
&+\kappa g(Q_{*,h}^{n+\frac{1}{2}}, s_{*}^{n+\frac{1}{2}})\sum_{i,j=1}^{3}\Big<Q_{h}^{i j, n+\frac{1}{2}}-Q_{*,h}^{ij,n+\frac{1}{2}},
\frac{Q_{h}^{i j, n+1}-Q_{h}^{i j, n}}{\tau}\Big>_{h},\label{Q43b}\\
\textcolor{black}{s^{n+1}=}&\textcolor{black}{\max\left\{\tilde{s}^{n+1},-C_\ast-E_{el}[Q^{n+1}_{h}]\right\}}, \label{Q43c}
\end{align}
\end{subequations}
where the operator $\mathcal{D}^{c}_{h}$ is defined in \eqref{optor} and $Q^{ij,n+1/2}_{h}:=(Q^{ij,n+1}_{h}+Q^{ij,n}_{h})/2$ for $i,j=1,2,3$. Here, $Q_{*,h}^{n+\frac{1}{2}}$ and $s_{*}^{n+\frac{1}{2}}$ are updated by the sESAV1 scheme \eqref{Q23}, that is $[Q_{*,h}^{n+\frac{1}{2}},s_{*}^{n+\frac{1}{2}}]=\mbox{sESAV1}(Q^{n}_{h},s^{n},\tau/2).$

Using a similar approach reported in \cite{gudibanda2022convergence}, we can deduce that the above scheme  also maintains the trace-free and symmetry properties of $Q$. Moreover, it is unconditionally stable in the sense of a modified energy, as stated in the following theorem. The proof is similar to that of Theorem \ref{thm2}, so we omit it here and leave it for interested readers.

\begin{thm}\label{thm4}
For any $\kappa\geq 0$, the \textnormal{sESAV2} scheme is unconditionally energy-dissipative in the sense that
$$E_h[Q_{h}^{n+1},s^{n+1}]\leq E_h[Q_{h}^n,s^n],\quad n\geq0.$$
\end{thm}

\textcolor{black}{
	\begin{remark}
		By defining the variable 
		$s^{n+1}$ in \eqref{Q23c}, we can readily establish that 
		\begin{equation}
			-C_*\leq s^{n+1}+E_{el}[Q^{n+1}_{h}]=E_h[Q_{h}^{n+1},s^{n+1}]\leq  E_h[Q^0_{h},\mathcal{E}_{1,h}[Q^{0}_{h}]],
		\end{equation}
		indicating that the total energy $E_h[Q_{h}^{n+1},s^{n+1}]$ is bounded for both sESAV1 and sESAV2 schemes.
	\end{remark}
}

\subsection{The discrete maximum bound principle}
As stated in Lemma \ref{lemma1}, the MBP property of $Q$ in Frobenius-norm holds for the Q-tensor model \eqref{Q6} in cases of $d=2$ and $d=3$ with $L_2+L_3=0$. For these two cases, the sESAV1 scheme \eqref{Q23} can be rewritten as follows: Given the initial values of $(Q_{h}^0, s^0)$, to find $(Q_{h}^{n+1}, s^{n+1})\in\mathbf{S}^{(d)}_{h}\times\mathbb{R},\;d=2,3, \; n\geq0,$ such that
\begin{subequations}\label{Q24}
\begin{align}
\textcolor{black}{\frac{Q_{h}^{i j, n+1}-Q_{h}^{i j, n}}{\tau}}
=& L \Delta_h Q_{h}^{i j, n+1}+g(Q_{h}^n, s^n)\big(f(Q_{h}^{ij,n})-\kappa(Q_{h}^{i j, n+1}-Q_{h}^{i j, n})\big),\label{Q24a}\\
\textcolor{black}{\tilde{s}}^{n+1}=&s^{n}-g(Q_{h}^n, s^n)\sum^3_{i,j=1}\Big<f(Q_{h}^{ij, n}), Q_{h}^{i j, n+1}-Q_{h}^{i j, n}
\Big>_{h},\label{Q24b}
\\
\textcolor{black}{s^{n+1}=}&\textcolor{black}{\max\left\{\tilde{s}^{n+1},-C_\ast-E_{el}[Q^{n+1}_{h}]\right\}}, \label{Q24c}
\end{align}
\end{subequations}
with $L=L_1+\frac{L_2+L_3}{2},$ where we have used the fact that, for the tensor-function $Q\in\mathcal{S}^{(2)}$, it holds
\begin{equation*}
\begin{aligned}
\sum_{k=1}^{2} \big(\partial_j \partial_k Q^{i k}
+\partial_i \partial_k Q^{j k}\big)
- \sum_{k,l=1}^{2}\partial_l \partial_k Q^{l k}\delta^{i j}=\Delta Q^{i,j}, \quad i,j=1,2.
\end{aligned}
\end{equation*}
For convenience, we denote the scheme \eqref{Q24} as $[Q^{n+1}_{h},s^{n+1}]=\mbox{MBP-sESAV1}(Q^{n}_{h},$ $s^{n},\tau).$
Furthermore, we can rewrite it in vector form for $d=2,3$, as follows:
\begin{subequations}\label{eq5}
\begin{align}		
\mathbf{G}_{h}\vQ_{h}^{i j, n+1}
=&\frac{1}{\tau}\vQ_{h}^{i j, n}
+g(Q_{h}^n, s^n)
\big[\kappa \vQ_{h}^{i j, n}+f(\vQ_{h}^{ij,n})\big],\label{eq5a}\\
\textcolor{black}{\tilde{s}}^{n+1}=&s^{n}-g(Q_{h}^n, s^n)\sum_{i,j=1}^{3}\Big<f(\vQ_{h}^{ij, n}), \vQ_{h}^{i j, n+1}-\vQ_{h}^{i j, n}
\Big>_{h},\label{eq5b}
\\
\textcolor{black}{s^{n+1}=}&\textcolor{black}{\max\left\{\tilde{s}^{n+1},-C_\ast-E_{el}[Q^{n+1}_{h}]\right\}}, \label{eq5c}
\end{align}
\end{subequations}
where $\mathbf{G}_{h}:=((1/\tau+\kappa g(Q_{h}^n, s^n))I
-L D_h
)$, and $\vec{Q}^{ij,n}_{h}\in{\mathbb R}^{(M+1)^d}$ is the vector representation of $Q^{ij,n}_{h}\in\mathbf{E}^{0}_{h}$, with the elements organized in the order of the $x_1,\cdots,x_d$ directions.
Here $f(\vQ_{h}^{ij, n})\in{\mathbb R}^{(M+1)^d}$ represents the element-wise function of $\vQ_{h}^{ij, n}$ and tensor form $D_{h}$ of the discrete $\Delta_{h}$ is given by
\bq\label{T_D}
D_{h}=
\begin{cases}
\begin{array}{r@{}l}
&I\otimes \Lambda_{h}+\Lambda_{h}\otimes I, \quad\qquad\qquad\qquad\qquad\quad~\mbox{ for } d=2;\\
&I\otimes I\otimes \Lambda_{h}+I\otimes \Lambda_{h}\otimes I+\Lambda_{h}\otimes I\otimes I,\quad \mbox{ for }d=3.\\
\end{array}
\end{cases}
\eq
with $I$ denoting the identity matrix (with the matched dimensions) and $\Lambda_{h}$ is the discrete matrix of the Laplace operator in 1D.
Moreover, we denote by $\vQ_{h}:=[\vQ^{11}_{h},\cdots,$ $\vQ^{1d}_{h};\cdots;\vQ^{d1}_{h},\cdots,\vQ^{dd}_{h}]$ the vector form of the grid tensor-function $Q_{h}\in\mathbf{S}^{(d)}_{h}$ for $d=2,3$.
Analogous to the Frobenius-norm and the $L^{\infty}$-norm for tensor function $Q\in L^p(\Omega \rightarrow \mathbb{R}^{d \times d})$, the discrete Frobenius-norm and $L^{\infty}$-Frobenius-norm	for the grid tensor-function $\vQ_{h}$ are given respectively by
\beq
\big|\vQ_{h}\big|^{2}_{F}=\sum_{i,j=1}^{d}\Big(\vQ^{ij}_{h}\Big)^{*2},\quad
\big\|\vQ_{h}\big\|_{F,\infty}=\Big\|\big|\vQ_{h}\big|_{F}\Big\|_{\infty},
\eeq
where $\Big(\vQ^{ij}_{h}\Big)^{*2}:=\vQ^{ij}_{h}.*\vQ^{ij}_{h}.$
Before investigating the MBP preservation of the proposed scheme \eqref{Q24}, we first review some useful lemmas.
\begin{lem}
[\textcolor{black}{\cite{CMPX19}, \cite{HWZ22}, Lemma 3.1 in \cite{STY16}}]\label{lemma2}
Let $B=(b_{kl})$ be a real $N\times N$ matrix and $A=\alpha I-B$ with $\alpha >0$. If $B$ is a negative diagonally dominant matrix, i.e.,
\begin{equation}
b_{kk}\leq 0\quad  \text{and}\quad b_{kk}+\sum_{l\neq k}|b_{kl}|\leq 0,\ \ 1\leq k\leq N,
\end{equation}
then, $A$ is invertible and satisfies
$\|A^{-1}\|_{\infty}\leq 1/\alpha.$
{Furthermore, if $b_{k,l}\geq0$ for $k\neq l$ and $1\leq k,l\leq N$, then it holds $A^{-1}\geq 0$, which means that all elements of $A^{-1}$ are non-negative.}
\end{lem}

\begin{lem}\label{lem2}
Assume $A:=\big(a_{kl}\geq0\big)$ is a real non-negative $N\times N$ matrix, and $\vPhi^{ij},\; i,j=1,\cdots,d$ are real $N\times1$ vectors. Then, it holds
\beq
\sqrt{\sum_{i,j=1}^{d}\big(A\vPhi^{ij}\big)^{*2}}\leq A \sqrt{\sum_{i,j=1}^{d}\big(\vPhi^{ij}\big)^{*2}}.
\eeq
\end{lem}

\begin{proof}
The well-known Minkowski inequality is stated as follows:
\beq
\Big(\sum_{k=1}^{m}|x_{k}+y_{k}|^{p}\Big)^{\frac{1}{p}}\leq \Big(\sum_{k=1}^{m}|x_{k}|^{p}\Big)^{\frac{1}{p}}+\Big(\sum_{k=1}^{m}|y_{k}|^{p}\Big)^{\frac{1}{p}},\ \ 1\leq p<\infty,
\eeq
for any real numbers $x_i,y_i,(i=1,\cdots,m).$
Then, using it with $p=2$ for $m_1$ times, we can deduce that
\bq \label{Minkowski_ADD}
\Big(\sum_{l=1}^{m_2}\Big(\sum_{k=1}^{m_1}\alpha_{k}^{l}\Big)^{2}\Big)^{\frac{1}{2}}
\leq\sum_{k=1}^{m_1}\Big(\sum_{l=1}^{m_2}(\alpha^{l}_{k})^{2}\Big)^{\frac{1}{2}}.
\eq
Let the $N\times 1$ vectors $\vec{b}:=\sqrt{\sum_{i,j=1}^{d}\big(A\vPhi^{ij}\big)^{*2}}$ and $\vec{c}:=\sqrt{\sum_{i,j=1}^{d}\big(\vPhi^{ij}\big)^{*2}}$. Then, using \eqref{Minkowski_ADD}, we deduce that
\beqal
\vec{b}_{k}=\sqrt{\sum_{i,j=1}^{d}\Big(\sum_{l=1}^{N}a_{kl}\vPhi^{ij}_{l}\Big)^{2}}\leq\sum_{l=1}^{N}\sqrt{\sum_{i,j=1}^{d}a^{2}_{k,l}\big(\vPhi^{ij}_{l}\big)^{2}}=\sum_{l=1}^{N}a_{k,l}\sqrt{\sum_{i,j=1}^{d}\big(\vPhi^{ij}_{l}\big)^{2}}=\vec{c}_{k},
\eeqal
for $1\leq k\leq N$, which completes the proof.
\end{proof}

\begin{lem}\label{lem1}
For any $Q\in\mathcal{S}^{(d)}, d=2,3,$ assume that $\|Q\|_{L^{\infty}(\Omega)}\leq \eta^{(d)}$ and  $\kappa\geq \max\{a+c\big(\eta^{(d)}\big)^{2},\; 0\}$ with $\eta^{(d)}$ defined in \eqref{equ1-1} and \eqref{equ1}, then it holds that
\bq \label{equ5}
\dps\sqrt{\sum_{i,j=1}^{d}\big(\kappa Q^{ij}+f(Q^{ij})\big)^{2}}\leq \kappa|Q|_{F}+\overline{f}(|Q|_{F}),
\eq
with the function $\overline{f}(\cdot)$ defined by
\bq\label{equ6}
\overline{f}(\xi)=-a \xi+\frac{b}{\sqrt{6}}\xi^2-c\xi^3,
\eq
where $b=0$ for $d=2$ as discussed in Remark \ref{remk1}.
\end{lem}
\begin{proof}
For the 2D case, we use  $\kappa\geq \max\{a+c\big(\eta^{(2)}\big)^{2},\; 0\}$, $\|Q\|_{L^{\infty}(\Omega)}\leq \eta^{(2)}$ and the definition of $f(\cdot)$ in \eqref{Q20} to obtain that
\beqal
\sqrt{\sum_{i,j=1}^{d}\big(\kappa Q^{ij}+f(Q^{ij})\big)^{2}}&=\sqrt{(\kappa -a-c\operatorname{tr}(Q^2))^2\sum_{i,j=1}^{d}(Q^{ij})^{2}}\\
&=(\kappa -a-c|Q|^{2}_{F})|Q|_{F}=\kappa|Q|_{F}+\overline{f}(|Q|_{F}).
\eeqal
Then we obtain the desired result for $d=2$.

For the 3D case, it holds
\bq\label{equ3}
\begin{aligned}
&\dps\sqrt{\sum_{i,j=1}^{d}\big(\kappa Q^{ij}+f(Q^{ij})\big)^{2}}\\
\leq&\dps\sqrt{(\kappa -a-c\operatorname{tr}(Q^2))^2\sum_{i,j=1}^{d}(Q^{ij})^{2}}+\sqrt{\sum_{i,j=1}^{d}b^{2}\big((Q^{2})^{ij}-\frac{1}{3}\mbox{tr}(Q)\delta^{ij}\big)^{2}}\\
=&\dps(\kappa -a-c|Q|^{2}_{F})|Q|_{F}+b\Big|Q^2 -\frac{1}{3}\operatorname{tr}(Q^2)I\Big|_{F},
\end{aligned}
\eq
where we have used $\kappa\geq \max\{a+c\big(\eta^{(3)}\big)^{2},\; 0\}$ and $\|Q\|_{L^{\infty}(\Omega)}\leq \eta^{(3)}$.
Since the tensor-function $Q\in\mathcal{S}^{(3)}$, then we can write $Q=P^TAP$ with $P^{T}P=I$ and $A=\mbox{diag}(\lambda_1, \lambda_2, \lambda_3)$ being a diagonal matrix with $\lambda_3=-\lambda_1-\lambda_2$.
Furthermore, it follows from
$Q^2=P^TA^2P$ and $\operatorname{tr}(Q^2)=\operatorname{tr}(A^2)$ that,
\bq\label{equ4}
\begin{aligned}
&|Q^2 -\frac{1}{3}\operatorname{tr}(Q^2)I|_{F}
=|A^2 -\frac{1}{3}\operatorname{tr}(A^2)I|_{F}\\
&\quad=\frac{1}{3}\sqrt{(2\lambda_1^2-\lambda_2^2-\lambda_3^2)^2
+(-\lambda_1^2+2\lambda_2^2-\lambda_3^2)^2
+(-\lambda_1^2-\lambda_2^2+2\lambda_3^2)^2}\\
&\quad=\frac{1}{\sqrt{6}}\big[\lambda_1^2+\lambda_2^2+(\lambda_1+\lambda_2)^2\big]=\frac{1}{\sqrt{6}}|Q|^{2}_{F},
\end{aligned}
\eq
where we have used
\beqal
\lambda_1^2\lambda_2^2+\lambda_1^2\lambda_3^2+\lambda_2^2\lambda_3^2=\frac{1}{2}\big[\lambda_1^4+\lambda_2^4+(\lambda_1+\lambda_2)^4\big]=\frac{1}{4}(\lambda_1^2+\lambda_2^2+(\lambda_1+\lambda_2)^2)^2.\\ 
\eeqal
Combining \eqref{equ3} and \eqref{equ4} gives  the desired result \eqref{equ5} for $d=3$.
\end{proof}
\begin{remark}
For the 2D case $(b=0)$, we have $\overline{f}(\xi)=(-a-c\xi^{2})\xi$ from \eqref{equ6}. Thus, it follows from $c>0$ and the definition of $\eta^{(2)}\geq \sqrt{\frac{a_{-}}{c}}$ in \eqref{equ1-1} that $\overline{f}(\eta^{(2)})\leq0$.
Similarly, we can deduce that $\overline{f}(\eta^{(3)})\leq0$ for the 3D case.
\end{remark}

\begin{lem}\label{lemma4}
Assume $\kappa\geq \dps\max_{\xi\in[0,\eta^{(d)}]}|\overline{f}'(\xi)|$ with $\eta^{(d)}$ and $\textcolor{black}{\overline{f}(\cdot)}$ defined in \eqref{equ1} and \eqref{equ6}, respectively. Then it holds that
\bq\label{equa1}
|\kappa\xi+\overline{f}(\xi)|\leq \kappa \eta^{(d)},\qquad\forall\xi\in[0,\eta^{(d)}].
\eq
\end{lem}

\begin{thm}\label{thm3}
For $\eta^{(d)}$ defined in \eqref{equ1-1} and \eqref{equ1}, assume that the stabilizing parameter $\kappa$ satisfies
\bq\label{equa5}
\dps\kappa\geq \max\big\{a+c\big(\eta^{(d)}\big)^{2},\; \max_{\xi\in[0,\eta^{(d)}]}|\overline{f}'(\xi)|\big\}.
\eq
Then, the \textnormal{MBP-sESAV1} scheme \eqref{Q24} unconditionally preserves the \textnormal{MBP} for $\{Q_{h}^n\}$:
\beqal
\|\vQ_{h}^n\|_{F,\infty}\leq\eta^{d},\ \forall n\geq0.
\eeqal
\end{thm}

\begin{proof}
It follows from the definition of $\eta^{(d)}$ in \eqref{equ1-1} and \eqref{equ1} that $\|\vQ^{0}_{h}\|_{F,\infty}\leq\eta^{(d)}.$
We suppose $\|\vQ_{h}^k\|_{F,\infty}\leq\eta^{(d)}$ for any $0\leq k\leq n$, which also implies that $0\leq|\vQ_{h}^{k}|_{F}\leq\eta^{(d)}$ for $0\leq k\leq n$. Next, we will show $\|Q^{n+1}_{h}\|_{F,\infty}\leq\eta^{(d)}.$
From \eqref{eq5}, we have
\begin{equation}\label{Q150}
\vQ_{h}^{i j, n+1}
=\mathbf{G}^{-1}_{h}\Big(\frac{1}{\tau}\vQ_{h}^{i j, n}+g(Q_{h}^n, s^n)\big[\kappa\vQ_{h}^{i j, n}+f(\vQ_{h}^{ij,n})\big]\Big),
\end{equation}
Furthermore, we use the above equality, the definition of $\mathbf{G}_{h}$, Lemmas \ref{lemma2}-\ref{lem1} to obtain
\begin{align*}
|\vQ^{n+1}_{h}|_{F}=&\sqrt{\sum_{i,j=1}^{d}\big(\vQ_{h}^{i j, n+1}\big)^{*2}}\\
=&\sqrt{\sum_{i,j=1}^{d}\Big[\mathbf{G}^{-1}_{h}\Big(\frac{1}{\tau}\vQ_{h}^{i j, n}+g(Q_{h}^n, s^n)\big[\kappa\vQ_{h}^{i j, n}+f(\vQ_{h}^{ij,n})\big]\Big)\Big]^{*2}}\\
\leq&\mathbf{G}^{-1}_{h}\left[\frac{1}{\tau}\sqrt{\sum_{i,j=1}^{d}\big(\vQ_{h}^{i j, n}\big)^{*2}}+g(Q_{h}^n, s^n)\sqrt{\sum_{i,j=1}^{d}\big[\kappa\vQ_{h}^{i j, n}+f(\vQ_{h}^{ij,n})\big]^{*2}}\right]\\
\leq&\mathbf{G}^{-1}_{h}\Big[\frac{1}{\tau}|\vQ^{n}_{h}|_{F}+g(Q_{h}^n, s^n)\big[\kappa|\vQ_{h}^{ n}|_{F}+\overline{f}(|\vQ_{h}^{n}|_{F})\big]\Big].
\end{align*}
Thus, together with Lemma \ref{lemma2} and Lemma \ref{lemma4}, we derive
\beqal
\||\vQ^{n+1}_{h}|_{F}\|_{\infty}\leq& \Big\|\mathbf{G}^{-1}_{h}\Big[\frac{1}{\tau}|\vQ^{n}_{h}|_{F}+g(Q_{h}^n, s^n)\big[\kappa|\vQ_{h}^{ n}|_{F}+\overline{f}(|\vQ_{h}^{n}|_{F})\big]\Big]\Big\|_{\infty}\\
\leq& \|\mathbf{G}^{-1}_{h}\|_{\infty}\Big[\frac{1}{\tau}\|\vQ^{n}_{h}\|_{F,\infty}+g(Q_{h}^n, s^n)\big\|\kappa|\vQ_{h}^{ n}|_{F}+\overline{f}(|\vQ_{h}^{n}|_{F})\big\|_{\infty}\Big]\\
\leq&\Big(\frac{1}{\tau}+\kappa g(Q_{h}^n, s^n)\Big)^{-1}
\Big(\frac{1}{\tau}+\kappa g(Q_{h}^n, s^n)\Big)\eta^{(d)}=
\eta^{(d)},
\eeqal
which completes the proof.
\end{proof}

\textcolor{black}{
Building upon the stability analysis presented in \cite{ju2022generalized}, we can first derive estimates for the discrete $H^1$ and $H^2$ norms of the numerical solution $Q^{n}_h$. Specifically, we have
\begin{lem}\label{lem4}
Suppose  $\kappa$ satisfies \eqref{equa5}, there exists a constant $M>0$ that depends on $C_1,|\Omega|, T$, $Q_0, \kappa, L$, and $\left\|\bm{f}\right\|_{C\{Q:
		|Q|_{F}\leq\eta^{(d)}\}}$, such that
	\begin{equation*}
	\left\|\nabla_h Q_h^{n+1}\right\|^2_h
		+\tau\sum_{k=0}^{n}\left(\left\| \frac{Q_h^{k+1}-Q_h^{k}}{\tau}\right\|_h^2
		+\left\|\Delta_h Q_h^{k+1}\right\|^2_h\right) \leq M, \quad 0 \leq n \leq\lfloor T / \tau\rfloor-1 .
	\end{equation*}
\end{lem}
}
\textcolor{black}{
\begin{cor}
If $\kappa$ satisfies the condition \eqref{equa5},
then there exists a positive constant $G^\ast$ such that $0<g(Q^k_{h}, s^k)<G^\ast$ and $s^{n+1}$ is bounded from below.
\end{cor}
\begin{proof}
 Using the definition of $g(\cdot,\cdot)$ in \eqref{Q23}, Theorem \ref{thm2} and \eqref{equa3},  we can derive that
 \bq\label{eqn1}
 g(Q^n, s^n) \leq \exp \{ E_h[Q^0_{h},\mathcal{E}_{1,h}[Q^{0}_{h}]]+C_*\}:=G^{*},\quad  n\leq 0.
 \eq
Based on \eqref{Q24b}, we can derive
\begin{align}\label{Q222}
\tilde{s}^{k+1}\geq \tilde{s}^k-\frac{\left(g\left(Q_h^k, s^k\right)\right)^2\tau}{2} \left\| \bm{f}\left(Q_h^k\right)\right\|_h^2- \frac{\tau}{2} \left\| \frac{Q_h^{k+1}-Q_h^{k}}{\tau}\right\|_h^2,
\end{align}
By summing up \eqref{Q222} over $k$ from $0$ to $n$, we obtain
\begin{align*}
	\tilde{s}^{n+1}\geq& s^0-\frac{(G^\ast)^2T}{2 } \left\|\bm{f}\right\|^2_{C\{Q:
		|Q|_F\leq\eta^{(d)}\}}|\Omega|- \frac{\tau}{2} \sum^{n}_{k=0} \left\| \frac{Q_h^{k+1}-Q_h^{k}}{\tau}\right\|_h^2\\
		\geq&-C_\ast
		-\frac{1+\kappa\tau G^\ast}{2}
		\left\|\nabla_h Q_h^0\right\|_h^2
		-\frac{L+2}{2L}(G^\ast)^2C_1T,
\end{align*}
when combined with \eqref{Q24c}, this concludes the proof.
\end{proof}
}

With a similar proving process of Corollary 3.5 in \cite{ju2022generalized}, we can derive the following corollary directly.

\begin{cor}\label{cor1}
	For fixed $h$ and $\eta^{(d)}$ defined in \eqref{equ1-1} and \eqref{equ1}, if the stabilizing parameter $\kappa$ satisfies the condition \eqref{equa5},
	then there exists 
	 $G^*>0$ such that $$G_{*}<g(Q^k_{h}, s^k) \leq G^*,\quad n\geq0.$$
\end{cor}

Next, we will investigate the MBP-preservation of the second-order in time scheme \eqref{Q43} for the cases of $d=2$ and $L_2+L_3=0$ for $d=3$, which can be written as follows:
given the initial values of ($Q_{h}^0$, $s^0$), we find $(Q_{h}^{n+1}\textcolor{black}{,} s^{n+1})\in \mathbf{S}^{d}_{h},\mathbb{R}, d=2,3,\; n\geq1$ such that
\begin{subequations}\label{Q46}
\begin{equation}\label{Q46a}
\begin{aligned}
	&	((\frac{2}{\tau}+\kappa g(Q_{*,h}^{n+\frac{1}{2}}, s_{*}^{n+\frac{1}{2}}))I
	-L \Delta_h
	)
	Q_{h}^{i j, n+1}=((\frac{2}{\tau}-\kappa g(Q_{*,h}^{n+\frac{1}{2}}, s_{*}^{n+\frac{1}{2}}))I\\
		&\quad+ L \Delta_h
	)Q_{h}^{i j, n}+2g(Q_{*,h}^{n+\frac{1}{2}}, s_{*}^{n+\frac{1}{2}})
	(\kappa Q_{*,h}^{ij,n+\frac{1}{2}}+f(Q_{*,h}^{ij,n+\frac{1}{2}})
	),
\end{aligned}
\end{equation}
\begin{equation}\label{Q46b}
\begin{aligned}
	\textcolor{black}{\tilde{s}}^{n+1}=&s^{n}-	g(Q_{*,h}^{n+\frac{1}{2}}, s_{*}^{n+\frac{1}{2}})\sum_{i,j=1}^{3}\big<f(Q_{*,h}^{ij,n+\frac{1}{2}}), Q_{h}^{i j, n+1}
	- Q_{h}^{i j, n}\big>_h \\
&+\kappa g(Q_{*,h}^{n+\frac{1}{2}}, s_{*}^{n+\frac{1}{2}})\sum_{i,j=1}^{3}\big<Q_{h}^{i j, n+\frac{1}{2}}-Q_{*,h}^{ij,n+\frac{1}{2}}, Q_{h}^{i j, n+1}
	- Q_{h}^{i j, n}\big>_h,\\
\end{aligned}
\end{equation}
\begin{equation}
\textcolor{black}{s^{n+1}=}\textcolor{black}{\max\left\{\tilde{s}^{n+1},-C_\ast-E_{el}[Q^{n+1}_{h}]\right\}},
\end{equation}
\end{subequations}
where $[Q_{*,h}^{n+\frac{1}{2}},s_{*}^{n+\frac{1}{2}}]=\mbox{MBP-sESAV1}(Q^{n}_{h},s^{n},\tau/2).$ For convenience, we will refer to the above scheme as the MBP-sESAV2 scheme hereafter. Then, the vector form of the  MBP-sESAV2 scheme \eqref{Q46} reads
\begin{subequations}\label{Q46-1}
\begin{align}		
\widehat{\mathbf{G}}_{h}\vQ_{h}^{i j, n+1}=&\overline{\mathbf{G}}_{h}\vQ_{h}^{i j, n}+2g(Q_{*,h}^{n+\frac{1}{2}}, s_{*}^{n+\frac{1}{2}})
(\kappa \vQ_{*,h}^{ij,n+\frac{1}{2}}+f(\vQ_{*,h}^{ij,n+\frac{1}{2}})
),\label{Q46-1a}\\
\textcolor{black}{\tilde{s}}^{n+1}=&s^{n}-g(Q_{*,h}^{n+\frac{1}{2}}, s_{*}^{n+\frac{1}{2}})\sum_{i,j=1}^{3}\big<f(\vQ_{*,h}^{ij,n+\frac{1}{2}}),\vQ_{h}^{i j, n+1}
-\vQ_{h}^{i j, n}\big>_h\label{Q46-1b}\\
 &+\kappa g(Q_{*,h}^{n+\frac{1}{2}}, s_{*}^{n+\frac{1}{2}})\sum_{i,j=1}^{3}\big<\vQ_{h}^{i j, n+\frac{1}{2}}-\vQ_{*,h}^{ij,n+\frac{1}{2}},  \vQ_{h}^{i j, n+1}
- \vQ_{h}^{i j, n}\big>_h,\notag\\
\textcolor{black}{s^{n+1}=}&\textcolor{black}{\max\left\{\tilde{s}^{n+1},-C_\ast-E_{el}[Q^{n+1}_{h}]\right\}},
\end{align}
\end{subequations}
with matrices $\widehat{\mathbf{G}}_{h}$ and $\overline{\mathbf{G}}_{h}$ defined respectively by
\bq\label{eqn2}
\widehat{\mathbf{G}}_{h}:=(\frac{2}{\tau}+\kappa g(Q_{*,h}^{n+\frac{1}{2}}, s_{*}^{n+\frac{1}{2}}))I-L D_h,\quad
\overline{\mathbf{G}}_{h}:=(\frac{2}{\tau}-\kappa g(Q_{*,h}^{n+\frac{1}{2}}, s_{*}^{n+\frac{1}{2}}))I+ L D_h,
\eq
where $[\vQ_{*,h}^{n+\frac{1}{2}},s_{*}^{n+\frac{1}{2}}]=\mbox{MBP-sESAV1}(\vQ^{n}_{h},s^{n},\tau/2).$

\begin{thm}\label{thm10}
For $\eta^{(d)}$ given in \eqref{equ1-1} and \eqref{equ1}, if $\kappa$ satisfies \eqref{equa5},
and
\begin{equation}\label{Q50}
\tau\leq (\frac{\kappa G^{\ast}}{2}+
\frac{2^{d-1}L}{h^2})^{-1},
\end{equation}
where $G^{\ast}>0$ is defined in Corollary \ref{cor1}.
Then the \textnormal{MBP-sESAV2} scheme \eqref{Q46}
preserves the \textnormal{MBP} for $\{Q_{h}^n\}$, i.e., $\|\vQ_{h}^n\|_{F,\infty}\leq\eta,$ is valid for $n\geq0$.
\end{thm}

\begin{proof}
Suppose $\|\vQ^{k}\|_{F,\infty}\leq \eta^{(d)}$ for $n\geq 0$, and $(\vQ_{*,h}^{n+\frac{1}{2}},s_{*}^{n+\frac{1}{2}})$ is given by the MBP-sESAV1 scheme \eqref{Q24}.
Then, we derive from Theorem \ref{thm3} and Corollary \ref{cor1} that $\|\vQ_{*,h}^{n+\frac{1}{2}}\|_{F,\infty}\leq\eta^{(d)}$ and $g(\vQ_{*,h}^{n+\frac{1}{2}},s_{*}^{n+\frac{1}{2}})\leq G^{*}$.
From the definitions of $D_h$ in \eqref{T_D} and $\widehat{G}_{h}$, it follows that
\beq
\widehat{G}^{-1}_{h}\geq0,\quad \|\widehat{G}^{-1}_{h}\|_{\infty}\leq (\frac{2}{\tau}+\kappa g(Q_{*,h}^{n+\frac{1}{2}}, s_{*}^{n+\frac{1}{2}}))^{-1},
\eeq
where we have used Lemma \ref{lemma2}. Moreover, under the condition \eqref{Q50} together with the definition of $D_{h}$, we derive that
\beq
\overline{\mathbf{G}}_{h}=(\frac{2}{\tau}-\kappa g(Q_{*,h}^{n+\frac{1}{2}}, s_{*}^{n+\frac{1}{2}}))I+ L D_h\geq0,
\eeq
which means all the elements of $\overline{\mathbf{G}}_{h}$ are non-negative. In addition, together with $\sum_{l=1}^{(M+1)^{d}}\big(D_{h}\big)_{kl}=0$ for any $k=1,2,\cdots,(M+1)^{d},$ we obtain
\beq
\|\overline{\mathbf{G}}_{h}\|_{\infty}=\frac{2}{\tau}-\kappa g(Q_{*,h}^{n+\frac{1}{2}}, s_{*}^{n+\frac{1}{2}}).
\eeq
Furthermore, using \eqref{Q46-1a}, Lemma \ref{lem2} and Lemma \ref{lem1} gives
\begin{equation}
\begin{aligned}		
	&|\vQ_{h}^{i j, n+1}|_{F} =\sqrt{\sum_{i,j}^{d}\big(\vQ_{h}^{i j, n+1}\big)^{*2}}\\
	=&\sqrt{\sum_{i,j}^{d}\Big[\widehat{\mathbf{G}}^{-1}_{h}\Big(\overline{\mathbf{G}}_{h}\vQ_{h}^{i j, n}
		+2g(Q_{*,h}^{n+\frac{1}{2}}, s_{*}^{n+\frac{1}{2}})
		\big(\kappa \vQ_{*,h}^{ij,n+\frac{1}{2}}+f(\vQ_{*,h}^{ij,n+\frac{1}{2}})
		\big)\Big)\Big]^{*2}}\\
	\leq&\widehat{\mathbf{G}}^{-1}_{h}\left[\overline{\mathbf{G}}_{h}\sqrt{\sum_{i,j}^{d}\big(\vQ_{h}^{i j, n}\big)^{*2}}
	+2g(Q_{*,h}^{n+\frac{1}{2}}, s_{*}^{n+\frac{1}{2}})
	\sqrt{\sum_{i,j}^{d}\big(\kappa \vQ_{*,h}^{ij,n+\frac{1}{2}}+f(\vQ_{*,h}^{ij,n+\frac{1}{2}})
		\big)^{*2}}\right]\\
	\leq&\widehat{\mathbf{G}}^{-1}_{h}\Big[\overline{\mathbf{G}}_{h}|\vQ_{h}^{i j, n}|_{F}
	+2g(Q_{*,h}^{n+\frac{1}{2}}, s_{*}^{n+\frac{1}{2}})
	\big(\kappa |\vQ_{*,h}^{ij,n+\frac{1}{2}}|_{F}+\overline{f}(|\vQ_{*,h}^{ij,n+\frac{1}{2}}|_{F})
	\big)\Big].
\end{aligned}
\end{equation}
Then, it follows from the above inequality, Lemma \ref{lemma2} and Lemma \ref{lemma4} that
\begin{equation*}
\begin{aligned}		
	&\big\||\vQ_{h}^{i j, n+1}|_{F}\big\|_{\infty}\\
	\leq&\Big\|\widehat{\mathbf{G}}^{-1}_{h}\Big[\overline{\mathbf{G}}_{h}|\vQ_{h}^{i j, n}|_{F}
	+2g(Q_{*,h}^{n+\frac{1}{2}}, s_{*}^{n+\frac{1}{2}})
	\big(\kappa |\vQ_{*,h}^{ij,n+\frac{1}{2}}|_{F}+\overline{f}(|\vQ_{*,h}^{ij,n+\frac{1}{2}}|_{F})
	\big)\Big]\Big\|_{\infty}\\
	\leq&\|\widehat{\mathbf{G}}^{-1}_{h}\|_{\infty}\Big[\|\overline{\mathbf{G}}_{h}\|_{\infty}\big\||\vQ_{h}^{i j, n}|_{F}\big\|_{\infty}
	+2g(Q_{*,h}^{n+\frac{1}{2}}, s_{*}^{n+\frac{1}{2}})
	\big\|\kappa |\vQ_{*,h}^{ij,n+\frac{1}{2}}|_{F}+\overline{f}(|\vQ_{*,h}^{ij,n+\frac{1}{2}}|_{F})
	\big\|_{\infty}\Big]\\
	\leq&\Big(\frac{2}{\tau}+\kappa g(Q_{*,h}^{n+\frac{1}{2}}, s_{*}^{n+\frac{1}{2}})\Big)^{-1}\Big[\Big(\frac{2}{\tau}-\kappa g(Q_{*,h}^{n+\frac{1}{2}}, s_{*}^{n+\frac{1}{2}})\Big)\eta^{(d)}
	+2g(Q_{*,h}^{n+\frac{1}{2}}, s_{*}^{n+\frac{1}{2}})\kappa\eta^{(d)}\Big]\\
	=&\eta^{(d)},
\end{aligned}
\end{equation*}
which completes the proof.
\end{proof}

We can  also derive that  $G_* \leq g(Q_{*,h}^{n+\frac{1}{2}}, s_{*}^{n+\frac{1}{2}}) \leq G^*$, and \textcolor{black}{  observe that $\left\|\nabla_h Q_h^{n+1}\right\|_h$ and $s^{n+1}$ are bounded}, which is similar to the analysis for the MBP-sESAV1 scheme \eqref{Q24} in Corollary \ref{cor1}. 
{\color{black}
\begin{remark}
The standard SAV method can be applied to the $Q$-tensor flow model \eqref{Q9} to derive a modified energy dissipative scheme. However, it appears challenging to preserve the MBP property using the standard SAV approach, as the positivity of the SAV term  $s^{n+1/2}/\sqrt{\mathcal{E}_{1,h}[Q^{n+1/2}_{*,h}]}$ cannot be readily guaranteed.
In contrast, the proposed sESAV method can unconditionally preserve a modified energy dissipation. Moreover, we can establish the MBP preservation of the proposed scheme due to the positivity of the well-designed ESAV function $g(Q^{n+1/2}_{*,h},s^{n+1/2}_{*}):=\exp(s^{n+1/2})/\exp(\mathcal{E}_{1,h}[Q^{n+1/2}_{*,h}])$. In the implementation, we update $Q^{n+1/2}_{*,h}$ and $s^{n+1/2}_{*}$ by using the first-order MBP-preserving sESAV scheme, denoted as  $[Q^{n+1/2}_{*,h},s^{n+1/2}_{*}]=\mbox{MBP-sESAV1}(Q^{n}_{h},s^{n},\tau/2)$.
Alternatively, one can define a different auxiliary variable such that the SAV function $g(Q,s)>0$ to obtain a modified energy dissipation and MBP-preserving scheme. A general definition of $g(\cdot,\cdot)$ has been provided in \cite{ju2022generalized} as follows 
\beq
g(Q,s)=\frac{\sigma(s)}{\sigma(\mathcal{E}[Q])},
\eeq
where $\sigma(\cdot)$  is a continuously differentiable, positive  function with $\sigma'(\cdot)\geq 0$  on $\mathbb{R}$.
\end{remark}
}

\section{Error Estimates}
In this section, we carry out rigorous error estimates for the fully discrete second-order MBP-sESAV scheme \eqref{Q46}.
For simplicity, we define the error functions as follows:
\begin{equation*}
\begin{aligned}
&e_Q^{ij, n}=Q_h^{ij, n}-Q^{ij}_e\left(t_n\right),
\quad
e_{*,Q}^{ij, n+\frac{1}{2}}=Q_{*,h}^{ij, n+\frac{1}{2}}-Q^{ij}_e\left(t_{n+\frac{1}{2}}\right),
\quad e_s^n=s^n-s_e\left(t_n\right),\\
&
\textcolor{black}{\tilde{e}_s^n=\tilde{s}^n-s_e\left(t_n\right)},
 \quad e_{*,s}^{n+\frac{1}{2}}=s_*^{n+\frac{1}{2}}-s_e\left(t_{n+\frac{1}{2}}\right),
\quad \textcolor{black}{\tilde{e}_{*,s}^{n+\frac{1}{2}}=\tilde{s}_*^{n+\frac{1}{2}}-s_e\left(t_{n+\frac{1}{2}}\right)},
\end{aligned}
\end{equation*}
where $Q^{ij}_e(t)$ and $s_e(t)$ represent the corresponding exact solutions, respectively.

Firstly, we shall provide estimates for $e_{*,Q}^{n+\frac{1}{2}}$ and $e_{*,s}^{n+\frac{1}{2}}$, which will be used in the derivation of error
estimate for the MBP-sESAV2 scheme \eqref{Q46}.

\begin{lem}\label{lemma42}
Under the assumption of Theorem \ref{thm3}, we have
\begin{equation}\label{Q888}
\|e_{*,Q}^{n+\frac{1}{2}}\|_h
\textcolor{black}{+\|\nabla_h e_{*,Q}^{n+\frac{1}{2}}\|_{h}}
+|e_{*,s}^{n+\frac{1}{2}}| \leq
C\big(\|e_Q^n\|_{h}+|e_s^n|+\tau(\tau+h^2)\big),
\end{equation}
where the constant $C$ depends on $C_*,|\Omega|,G^{*}, Q_e,$ and $\left\|f\right\|_{C^1\{Q:
|Q|_{\textcolor{black}{F}}\leq\eta^{\textcolor{black}{(d)}}\}}$.
\end{lem}

\begin{proof}
From \eqref{Q20} and \eqref{Q24}, we deduce  error equations with respect to
$e_{*,Q}^{ij,n+\frac{1}{2}}$
and $e_{*,s}^{n+\frac{1}{2}}$ as follows:

\begin{subequations}\label{Q87}
	\begin{align}
		\big(\frac{2}{\tau}+\kappa g\left(Q_{h}^n, s^n\right)\big)&(e_{*,Q}^{ij,n+\frac{1}{2}}-e_Q^{ij, n})-L\Delta_h e_{*,Q}^{ij,n+\frac{1}{2}}=J^{ij,n}_{1,h}+J^{ij,n}_{2,h}-R_{1, Q}^{ij,n}, \label{Q87a}\\
		\frac{\textcolor{black}{\tilde{e}}_{*,s}^{n+\frac{1}{2}}-e_s^{n}}{\tau/2}=&\dps-g(Q_{h}^n, s^n)\sum_{i,j=1}^3\Big< f(Q_{h}^{ij,n}), \frac{e_{*,Q}^{n+\frac{1}{2}}-e_Q^{ij, n}}{\tau/2} \Big>_h\label{Q87b}\\
&-\sum_{i,j=1}^3\Big<J^{ij,n}_{1,h}, \frac{Q^{ij}_e(t_{n+\frac{1}{2}})-Q^{ij}_e(t_n)}{\tau/2}\Big>_h-R_{1,s}^n ,\notag
\end{align}
\end{subequations}
where the grid tensor functions $J^{n}_{1,h}$ and $J^{n}_{2,h}$ for $1\leq i,j\leq d$ are given by
\bry
 J^{ij,n}_{1,h}&:=g\left(Q_{h}^n, s^n\right) f(Q_{h}^{ij, n})- f(Q^{ij}_e(t_{n+\frac{1}{2}})),\\
 J^{ij,n}_{2,h}&:=-\kappa g\left(Q_{h}^n, s^n\right)\big(Q^{ij}_e(t_{n+\frac{1}{2}})-Q^{ij}_e\left(t_n\right)\big),
 \ery
and the corresponding truncation errors $R_{1, Q}^{n}$ and $R_{1, s}^n$  {\color{black} with respect to the time derivative and Laplace operator discretizations of $Q$, and the time derivative discretization of $s$, respectively,}  can be estimated as follows:
\begin{equation}\label{Q88}
\|R_{1, Q}^n\|_h \leq C(\tau+h^2), \quad|R_{1, s}^n| \leq C\tau.
\end{equation}	
Taking the discrete inner product of \eqref{Q87a} with $e_{*,Q}^{ij,n+\frac{1}{2}}\textcolor{black}{-e_Q^{ij,n}}$, summing up $i,\; j$ from 1 to d, and using Young's inequality, we can obtain that
\textcolor{black}{
\begin{equation}\label{QQ222}
	\begin{aligned}
		&
		\frac{L}{2}\Big(\|\nabla_h e_{*,Q}^{n+\frac{1}{2}}\|_{h}^2-\|\nabla_h e_{Q}^{n}\|_{h}^2
	+\|\nabla_h e_{*,Q}^{n+\frac{1}{2}}-\nabla_h e_{Q}^{n}\|_{h}^2\Big)\\
&+	\Big(\frac{2}{\tau}+\kappa g(Q_{h}^n, s^n)\Big)
	\|e_{*,Q}^{n+\frac{1}{2}}-e_Q^{n}\|_{h}^2\\
		=&\sum_{i,j=1}^{3}\Big<J^{ij,n}_{1,h}+J^{ij,n}_{2,h}-R_{1*, Q}^{ij,n}, e_{*,Q}^{ij,n+\frac{1}{2}}-e_{Q}^{ij,n}\Big>_h \\
		\leq&\frac{\textcolor{black}{3}\tau}{2}\big(\|J^{n}_{1,h}\|^{2}_{h}+\|J^{n}_{2,h}\|^{2}_{h}+\|R_{1, Q}^{n}\|^{2}_{h}\big)+\frac{1}{2\tau}\| e_{*,Q}^{n+\frac{1}{2}}-e_{Q}^{n}\|^{2}_{h}.
	\end{aligned}
\end{equation}
}
\textcolor{black}{
Using the following identity
$(a-b)^2=a^2-b^2-2(a-b,b),$
we derive
	\begin{equation}
		\begin{aligned}
			&
			\frac{L}{2}\Big(\|\nabla_h e_{*,Q}^{n+\frac{1}{2}}\|_{h}^2-\|\nabla_h e_{Q}^{n}\|_{h}^2
			\Big)
			+	\kappa G_*
			(\|e_{*,Q}^{n+\frac{1}{2}}\|_h^2-\|e_Q^{n}\|_{h}^2)
				+	\frac{1}{\tau}
			\|e_{*,Q}^{n+\frac{1}{2}}-e_Q^{n}\|_{h}^2\\
		\leq&\frac{\textcolor{black}{3}\tau}{2}\big(\|J^{n}_{1,h}\|^{2}_{h}+\|J^{n}_{2,h}\|^{2}_{h}+\|R_{1, Q}^{n}\|^{2}_{h}\big)+
		\frac{\kappa^2 G^2_*\tau}{2}
	\|e_Q^{n}\|^2_h.
		\end{aligned}
	\end{equation}
}
Thus, together with \eqref{Q88} and Corollary \ref{cor1}, we deduce  that
\bq\label{Q889}
\|e_{*,Q}^{n+\frac{1}{2}}\|_{h}^2
\textcolor{black}{+\|\nabla_h e_{*,Q}^{n+\frac{1}{2}}\|_{h}^2}
\leq C\|e_Q^{n}\|_{h}^2+\tau^{2}\|J^{n}_{1,h}\|^{2}_{h}+C_{Q_{e}}\kappa G^{*}\tau^{4}+C\tau^{2}(\tau+h^{2})^{2}.
\eq
Noting that $\|Q_e(t_n)\|_{\infty}\leq \eta^{(d)},\; \|\vQ_h^{n}(t_n)\|_{F,\infty} \leq \eta^{(d)}$  and the uniform bounds of $s_e$ and $s^{n}$, we next give an estimation for the term $J^{n}_{1,h}$.
Denote the term $g(Q_{h}^n, s^n)-g(Q_e(t_n), s_e(t_n))$ by $J^{n}_{3,h}$, then we have
\begin{equation}
\begin{aligned}
	\big|J^{n}_{3,h}\big|&\leq\left|g\left(Q_{h}^n, s^n\right)-g\left(Q_{h}^n, s_e\left(t_n\right)\right)\right|+\left|g\left(Q_{h}^n, s_e\left(t_n\right)\right)-g\left(Q_e\left(t_n\right), s_e\left(t_n\right)\right)\right|\\
&:=A_1+A_2.
\end{aligned}
\end{equation}
In addition, we deduce that
\begin{equation*}
\begin{aligned}
	A_1&=\frac{1}{\exp \left\{\mathcal{E}_{1 h}\left(Q_{h}^n\right)\right\}}\left|\exp \left(s^n\right)-\exp \left\{s_e\left(t_n\right)\right\}\right|
	 \leq C_g\left|e_s^n\right|,\\
A_2	&	=\exp \left\{s_e\left(t_n\right)\right\}\left|\frac{1}{\exp \left\{\mathcal{E}_{1 h}\left(Q_{h}^n\right)\right\}}-\frac{1}{\exp \left\{\mathcal{E}_{1 h}\left(Q_e\left(t_n\right)\right)\right\}}\right| \\
	&	 \leq\exp \left\{\mathcal{E}\left(Q_h^0\right)+C_*\right\}
	\left|\mathcal{E}_{1 h}\left(Q_{h}^n\right)-\mathcal{E}_{1 h}\left(Q_e\left(t_n\right)\right)\right| 
		 \leq C_g\left\|e_Q^n\right\|_h,
\end{aligned}
\end{equation*}
where the constant $C_g$ depends on $C_*,\ |\Omega|.$
Hence, $J^{n}_{3,h}$ can be bounded by
\begin{equation}\label{Q76}
\big|J^{n}_{3,h}\big| \leq C_g(\|e_Q^n\|_h+|e_s^n|).
\end{equation}
Furthermore, we use triangle inequality and H\"older's inequality to obtain
\begin{equation}
\begin{aligned}\label{Q77}
\|J^{n}_{1,h}\|_{h}
	\leq&\|g(Q_{h}^n, s^n)\big(f(Q_{h}^n)- f(Q_e(t_n))\big)\|_h +\|J^{n}_{3,h} f(Q_e(t_n))\|_h\\
&+\| f(Q_e(t_n))- f(Q_e(t_{n+\frac{1}{2}}))\|_h\\
\leq&\|f^{\prime}\|_{ C\{Q:
		|Q|_{F}\leq\eta^{(d)}\}}\big[G^*\|e_Q^n\|_h+C_{Q_e}\tau\big]+\|f(Q_e(t_n))\|_h|J^{n}_{3,h}||\Omega|\\
	\leq& C(\|e_Q^n\|_h+|e_s^n|+\tau),
\end{aligned}
\end{equation}		
where the constant $C$ depends on $C_g,\; G^{*}$ and $\left\|f\right\|_{C^1\{Q:
|Q|\leq\eta\}}$.
Substituting \eqref{Q77} into \eqref{Q889}, we get
\begin{equation}\label{Q89}
\begin{aligned}
	\big\|e_{*,Q}^{n+\frac{1}{2}}\big\|_{h}^2
	\textcolor{black}{+\|\nabla_h e_{*,Q}^{n+\frac{1}{2}}\|_{h}^2}
	\leq
	C\left(\left\|e_Q^n\right\|_{h}^2+\tau^{2}\left|e_s^n\right|^2+\tau^{2}\left(\tau+h^2\right)^2\right).
\end{aligned}
\end{equation}
By multiplying \eqref{Q87b} by $\tau e_{*,s}^{n+\frac{1}{2}}$ and utilizing \eqref{Q88}, \eqref{Q77}, \eqref{Q89}, and Young's inequality, we can obtain that
\begin{equation*}
\begin{aligned}
	&\big(\big|\textcolor{black}{\tilde{e}}_{*,s}^{n+\frac{1}{2}}\big|^2
	-|e_s^n|^2
	+\big|\textcolor{black}{\tilde{e}}_{*,s}^{n+\frac{1}{2}}-e_s^n\big|^2\big) \\
	&~~=-\textcolor{black}{\tilde{e}}_{*,s}^{n+\frac{1}{2}}\Big[2 g(Q_{h}^n, s^n)\sum_{i,j=1}^{3}\Big < f(Q_{h}^{ij,n}), e_{*,Q}^{ij,n+\frac{1}{2}}-e_Q^{ij,n}\Big>_h\\
&~~\quad+2\sum_{i,j=1}^{3}\Big <J^{ij,n}_{1,h},Q^{ij}_e(t_{n+\frac{1}{2}})-Q^{ij}_e(t_n) \Big>_h
	+\tau R_{1,s}^n \Big]\\
	&~~\leq\frac{1}{2}|\textcolor{black}{\tilde{e}}_{*,s}^{n+\frac{1}{2}}|^2
	+C\Big[\|e_{*,Q}^{ij,n+\frac{1}{2}}\|^{2}_{h}+\|e_Q^{ij,n}\|^{2}_{h}
	+\|J^{n}_{1,h}\|^{2}_{h}\tau^{2}+ \tau^{2}(\tau+h^2)^2\Big]\\
&~~\leq\frac{1}{2}|\textcolor{black}{\tilde{e}}_{*,s}^{n+\frac{1}{2}}|^2
	+C\Big[\left\|e_Q^n\right\|_{h}^2+\tau^{2}\left|e_s^n\right|^2+\tau^{2}\left(\tau+h^2\right)^2\Big].
\end{aligned}
\end{equation*}
Thus, we derive that
\begin{equation}\label{Q90}
\begin{aligned}
	\big|\textcolor{black}{\tilde{e}}_{*,s}^{n+\frac{1}{2}}\big|^2\leq
	2\left|e_s^n\right|^2+C\Big[\left\|e_Q^n\right\|_{h}^2+\tau^{2}\left|e_s^n\right|^2+\tau^{2}\left(\tau+h^2\right)^2\Big].
\end{aligned}
\end{equation}
\textcolor{black}{If $\tilde{s}_{*}^{n+\frac{1}{2}}\geq-C_\ast-E_{el}[Q_{*,h}^{n+\frac{1}{2}}]$, then $s_{*}^{n+\frac{1}{2}}=\tilde{s}_{*}^{n+\frac{1}{2}}$, and thus $e_{*,s}^{n+\frac{1}{2}}=\tilde{e}_{*,s}^{n+\frac{1}{2}}$.
Conversely, if $\tilde{s}_{*}^{n+\frac{1}{2}} < -C_\ast-E_{el}[Q_{*,h}^{n+\frac{1}{2}}]$, we find that
	\begin{equation}\label{Q99}
	\tilde{s}_{*}^{n+\frac{1}{2}}+E_{el}[Q_{*,h}^{n+\frac{1}{2}}]\leq s_{*}^{n+\frac{1}{2}}+E_{el}[Q_{*,h}^{n+\frac{1}{2}}]=-C_\ast\leq s_e(t_{n+\frac{1}{2}})+E_{el}[Q_e(t_{n+\frac{1}{2}})],
	\end{equation}
which leads to
	\begin{equation}
		\tilde{s}_{*}^{n+\frac{1}{2}}\leq s_{*}^{n+\frac{1}{2}}\leq s_e(t_{n+\frac{1}{2}})+E_{el}[Q_e(t_{n+\frac{1}{2}})]-E_{el}[Q_{*,h}^{n+\frac{1}{2}}].
		\end{equation}
By incorporating \eqref{Q90}, we can obtain
\begin{equation}\label{Q999}
	\begin{aligned}
\big|e_{*,s}^{n+\frac{1}{2}}\big|^2\leq&
\big|\tilde{e}_{*,s}^{n+\frac{1}{2}}\big|^2+
\frac{L}{2}(\|\nabla Q_e(t_{n+\frac{1}{2}})\|^2-\|\nabla_h Q^{n+\frac{1}{2}}_{*,h}\|^2_{h})^2\\
\leq&\big|\tilde{e}_{*,s}^{n+\frac{1}{2}}\big|^2+C(\|\nabla_h e_{*,Q}^{k+\frac{1}{2}}\|_{h}^2+h^2)\\
\leq&2\left|e_s^n\right|^2+C\Big[\left\|e_Q^n\right\|_{h}^2+\tau^{2}\left|e_s^n\right|^2+\tau^{2}\left(\tau+h^2\right)^2\Big].
	\end{aligned}
\end{equation}
}
Combing \eqref{Q89} \textcolor{black}{and \eqref{Q999}} together, we obtain the desired results \eqref{Q888}.
\end{proof}

\begin{thm}
Under the assumption of Theorem \ref{thm10}, it holds for  \textnormal{MBP-sESAV\textcolor{black}{p}}  scheme \eqref{Q24} \textcolor{black}{for $p=1$} and \eqref{Q46} \textcolor{black}{for $p=2$} that
\begin{equation}\label{est4}
\left\|e_Q^n\right\|+\left\|\nabla_h e_Q^n\right\|+\left|e_s^n\right| \leq C\left(\tau^{\textcolor{black}{p}}+h^2\right), \quad 0 \leq n \leq\lfloor T / \tau\rfloor,
\end{equation}
where the constant $C$ depends on $C_*,|\Omega|,G_{*},G^{*}, Q_e,T$, and $\left\|f\right\|_{C^1\{Q:
	|Q|_{\textcolor{black}{F}}\leq\eta^{\textcolor{black}{(d)}}\}}$.
\end{thm}

\begin{proof}
	\textcolor{black}{
We only provide estimates for MBP-sESAV2 scheme, while omitting the proof for MBP-sESAV1 scheme \eqref{Q24} since their proofs are essentially similar.
}

It follows from \eqref{Q20} and \eqref{Q46} that the error equations with respect to
$e_Q^{ij,n+1}$and $e_s^{n+1}$ read as
\begin{subequations}\label{Q92}
\begin{equation}\label{Q92a}
	\begin{aligned}
		&\frac{e_Q^{ij, n+1}-e_Q^{ij, n}}{\tau}-L\Delta_h\frac{ e_Q^{ij, n+1}+e_Q^{ij, n}}{2}+\kappa g\big(Q_{*,h}^{n+\frac{1}{2}}, s_{*}^{n+\frac{1}{2}}\big)\frac{ e_Q^{ij, n+1}+e_Q^{ij, n}}{2}\\
		&~~=J^{ij,n}_{1*,h}+J^{ij,n}_{2*,h}+\kappa g\big(Q_{*,h}^{n+\frac{1}{2}}, s_{*}^{n+\frac{1}{2}}\big)e_{*,Q}^{ij, n+\frac{1}{2}} -R_{2, Q}^{ij,n},
	\end{aligned}
\end{equation}	
	\begin{equation}\label{Q92b}
	\begin{aligned}
		\frac{\textcolor{black}{\tilde{e}}_s^{n+1}-e_s^{n}}{\tau} =&-g(Q_{*,h}^{n+\frac{1}{2}}, s_{*}^{n+\frac{1}{2}})\sum_{i,j=1}^{3}\Big < f(Q_{*,h}^{ij,{n+\frac{1}{2}}}), \frac{e_Q^{ij, n+1}-e_Q^{ij, n}}{\tau} \Big>_h\\
		&-\sum_{i,j=1}^{3}\Big < J^{ij,n}_{1*,h}, \frac{Q^{ij}_e(t_{n+1})-Q^{ij}_e(t_n)}{\tau}\Big>_h -R_{2 ,s}^n,
	\end{aligned}
		\end{equation}
\end{subequations}
where the grid tensor functions $J^{n}_{1*,h}$ and $J^{n}_{2*,h}$ for $1\leq i,j\leq d$ are given by
\bry
 J^{ij,n}_{1*,h}&:=g(Q_{*,h}^{n+\frac{1}{2}}, s_{*}^{n+\frac{1}{2}}) f(Q_{*,h}^{ij, n+\frac{1}{2}})- f(Q^{ij}_e(t_{n+\frac{1}{2}})),\\
 J^{ij,n}_{2*,h}&:=-\kappa g\big(Q_{*h}^{n+\frac{1}{2}}, s^{n+\frac{1}{2}}_{*}\big)
		\Big(\frac{Q^{ij}_e(t_{n+1})+Q^{ij}_e(t_{n})}{2}-Q^{ij}_e(t_{n+\frac{1}{2}})\Big),
 \ery
and the truncation errors $R_{2 Q}^n$ and $R_{2 s}^n$ {\color{black} with respect to the time derivative and Laplace operator discretizations of $Q$, and the time derivative discretization of $s$, respectively} satisfy
\beq
\|R_{2,Q}^n\|\textcolor{black}{_h} \leq C(\tau^2+h^2), \quad|R_{2,s}^n| \leq C\tau^2.
\eeq
Taking the discrete inner product of \eqref{Q92a} with $e_Q^{ij,n+1}-e_Q^{ij,n}$ and summing up for $1\leq i,j \leq 3$, we arrive at
\begin{equation}\label{Q93}
\begin{aligned}
	&\frac{	1}{\tau}\|e_Q^{n+1}-e_Q^{ n}\|_{h}^2
	+\frac{1}{2}\big[\kappa G_{*}(\|e_Q^{n+1}\|_{h}^2-\|e_Q^n\|_{h}^2)+L(\|\nabla_h e_Q^{n+1}\|_{h}^2-\|\nabla_h e_Q^n\|_{h}^2)\big]\\
	&~~=\sum_{i,j=1}^{3}\Big <J^{ij,n}_{1*,h}+J^{ij,n}_{2*,h}+\kappa g\big(Q_{*,h}^{n+\frac{1}{2}}, s_{*}^{n+\frac{1}{2}}\big)e_{*,Q}^{ij, n+\frac{1}{2}} -R_{2, Q}^{ij,n}, e_Q^{ij,n+1}-e_Q^{ij,n}\Big>_h \\
	&~~\leq\frac{1}{2\tau}\|e_Q^{n+1}-e_Q^n\|_{h}^2
	+C\tau\big(\|J^{n}_{1*,h}\|^{2}_{h}+\|e_{*,Q}^{n+\frac{1}{2}}\|_{h}^2+(\tau^2+h^2)^2\big).
\end{aligned}
\end{equation}	
{where we have used Young's inequality and $0<G_{*}\leq g(Q_{*,h}^{n+\frac{1}{2}}, s_{*}^{n+\frac{1}{2}})\leq G^*$.
By following a similar process of deriving the estimate for $J^{n}_{1*,h}$ in Lemma \ref{lemma42}, we can obtain a similar estimate for $J^{n}_{1*,h}$ as follows:} 
\beq
\|J^{n}_{1*,h}\|^{2}_{h}\leq C(\|e_{*,Q}^{n+\frac{1}{2}}\|^{2}_h+|e_{*,s}^{n+\frac{1}{2}}|^{2}).
\eeq	
Together with \eqref{Q93}, we obtain
\begin{equation}\label{est1}
\begin{aligned}
	&\dps\frac{	1}{2\tau}\|e_Q^{n+1}-e_Q^{ n}\|_{h}^2+\frac{1}{2}\big[\kappa G_{*}(\|e_Q^{n+1}\|_{h}^2-\|e_Q^n\|_{h}^2)+L(\|\nabla_h e_Q^{n+1}\|_{h}^2-\|\nabla_h e_Q^n\|_{h}^2)\big]\\
&~~\leq C\tau\big(\|e_{*,Q}^{n+\frac{1}{2}}\|_{h}^2+|e_{*,s}^{n+\frac{1}{2}}|^{2}+(\tau^2+h^2)^2\big).\\
\end{aligned}
\end{equation}	
Multiplying \eqref{Q92b} by $2 \tau e_s^{n+1}$ and using the H\"older's inequality, yields
\brr\label{est2}
	&| \textcolor{black}{\tilde{e}}_s^{n+1}|^2
	-|e_s^n|^2
	+|\textcolor{black}{\tilde{e}}_s^{n+1}-e_s^n|^2\\
	&~~= \dps-2\textcolor{black}{\tilde{e}}_s^{n+1}\big[g(Q_{*,h}^{n+\frac{1}{2}}, s_{*}^{n+\frac{1}{2}})\sum_{i,j=1}^{3}\big < f(Q_{*,h}^{ij,{n+\frac{1}{2}}}), e_Q^{ij, n+1}-e_Q^{ij, n} \big>_h\\
		&~~\quad\dps-\sum_{i,j=1}^{3}\big < J^{ij,n}_{1*,h}, Q^{ij}_e(t_{n+1})-Q^{ij}_e(t_n)\big>_h -\tau R_{2 ,s}^n \big]\\
&~~\leq C\big[|\textcolor{black}{\tilde{e}}_s^{n+1}|\big\| e_Q^{ n+1}-e_Q^{n}\big\|_{h}+|\textcolor{black}{\tilde{e}}_s^{n+1}|\big\|J^{n}_{1*,h}\big\|_{h}\tau+\tau |\textcolor{black}{\tilde{e}}_s^{n+1}||R_{2 ,s}^n|\big] \\
&~~\dps\leq\frac{\| e_Q^{n+1}-e_Q^{ n}\|^{2}_{h}}{2\tau}+C\tau\big[|\textcolor{black}{\tilde{e}}_s^{n+1}|^{2}+\|e_{*,Q}^{n+\frac{1}{2}}\|^{2}_h+|e_{*,s}^{n+\frac{1}{2}}|^{2}+\tau^{2}(\tau+h^{2})^{2}\big].
\err
Summing up the inequalities \eqref{est1} and \eqref{est2}, and using the Lemma \ref{lemma42}, we derive
\begin{equation}\label{est3}
\begin{aligned}
	&\dps\frac{1}{2}\big[\kappa G_{*}(\|e_Q^{n+1}\|_{h}^2-\|e_Q^n\|_{h}^2)+L(\|\nabla_h e_Q^{n+1}\|_{h}^2-\|\nabla_h e_Q^n\|_{h}^2)\big]+| \textcolor{black}{\tilde{e}}_s^{n+1}|^2
	-|e_s^n|^2\\
&~~\leq C\tau\big(|\textcolor{black}{\tilde{e}}_s^{n+1}|^{2}+\|e_{*,Q}^{n+\frac{1}{2}}\|_{h}^2+|e_{*,s}^{n+\frac{1}{2}}|^{2}+(\tau^2+h^2)^2\big)\\
&~~\leq C\tau\big(|\textcolor{black}{\tilde{e}}_s^{n+1}|^{2}+\|e_Q^n\|_{h}^2+|e_s^n|^2+(\tau^2+h^2)^2\big).\\
\end{aligned}
\end{equation}	
\textcolor{black}{Similar to \eqref{Q99}-\eqref{Q999} in Lemma \ref{lemma42}, }
then \eqref{est4} can be obtained by using the discrete Gronwall's inequality.
\end{proof}	
\textcolor{black}{
	\begin{remark}
		It is noted that without the MBP property, we can still establish error results similar to those presented in \eqref{est4} by using a mathematical induction.
		For instance, considering the MBP-sESAV2 scheme \eqref{Q46} and assuming that 
		$\| e_Q^k\|_{H_h^2}\leq \tau+h$ for $1\leq k\leq n$, we can demonstrate  $\| e_Q^{n+1}\|_{H_h^2}\leq \tau+h$.
		By taking the discrete inner product of \eqref{Q87a} with $-(\Delta_h e_{*,Q}^{ij,n+\frac{1}{2}}-\Delta_h e_Q^{ij,n} )$ and  using similar estimates of \eqref{QQ222}-\eqref{Q999}, we can derive the $L^{\infty}$ boundedness of $ e_{*,Q}^{n+1/2} $ under the condition that $\tau$ and $h$ are sufficiently small,
		implying $0< G_*\leq g(Q_{*,h}^{n+1/2}, s_*^{n+1/2} ) \leq G^*$.
		Subsequently, by taking the discrete inner product of \eqref{Q92a} with $-(\Delta_h e_{Q}^{ij,n+1}-\Delta_h e_Q^{ij,n} )$, we readily achieve $\| e_Q^{n+1}\|_{H_h^2}\leq C(\tau^2+h^2)$. Correspondingly, it can be demonstrated that the inequality $\| e_Q^{n+1}\|_{H_h^2}\leq \tau+h$ is valid under the condition that $\tau$ and $h$ are sufficiently small, which leads to the $L^\infty$ boundedness of $Q^{n+1}$.
	\end{remark}
}

\section{Numerical Simulations}
In this section, we present various numerical experiments in 2D and 3D cases to verify the theoretical results derived in
the previous section and demonstrate the supremum norms, the energy stability and the accuracy of the proposed numerical schemes.
Moreover, the first-order scheme \eqref{Q23} and the second-order scheme \eqref{Q43} can be implemented by fast Discrete Sine Transform (DST) efficiently. For more information, see, e.g., \cite{ju2015fast}.
Unless otherwise specified, we choose  the second-order scheme \eqref{Q43} \textcolor{black}{with $C_*=1$}  in the following simulation.
\subsection{Convergence tests}
In this subsection, we verify the accuracy of the proposed  
 \mbox{MBP-sESAV2} schemes.
In the following simulations, the constant $\kappa$ is determined by $\kappa=2$, and we take the model parameter values and the initial condition as
\bry
	       &L_d=1,\quad L=1.0\times 10^{-3},\quad a=-0.25,\quad b=1, \quad c=1,\\
		&Q_0=\textbf{n}_0\textbf{n}_0^T-\frac{\|\textbf{n}_0\|_{h}^2}{2}\textbf{I},\quad
		\textbf{n}_0=(\sin(2\pi x)\sin(2\pi y),\
		\sin(2\pi x)\sin(2\pi y))^T.
\ery

We first test convergence rates of spatial discretization for \eqref{Q43}.
Particularly, $\|e_\xi\|=\|\xi_h-\xi_{h/2}\|$ determines the errors between the two distinct grid spacings $h$ and $h/2$.
The reference solution is \textcolor{black}{obtained using a very fine temporal discretization} with $2\times 10^4$ time steps.
The errors at $T=1$, measured in the $L^2$ norm with \textcolor{black}{various spatial step sizes} are shown in Table \ref{table2}, where we can get the second-order convergence in the \textcolor{black}{spatial} direction clearly.

Next, we give convergence rates of \textcolor{black}{temporal discretization} for \eqref{Q43}.
The reference solution is computed \textcolor{black}{using a fine spatial discretization with} $256$ grid points in each direction.
The numerical results for the MBP-sESAV2 schemes are reported in Table \ref{table4}, which give solid supporting evidence for the expected second-order convergence in time.
\textcolor{black}{To save space, we omit the details for the first-order scheme, as it also exhibits the desired convergence rates.}

\begin{table}[!t]
	\renewcommand{\arraystretch}{1.1}
	\small
	\centering
	\caption{Error and convergence rates \textcolor{black}{of spatial discretization} for \textcolor{black}{\textnormal{MBP-sESAV2}} scheme \eqref{Q43}.
	}\label{table2}
	\scalebox{0.86}{\begin{tabular}{p{1.0cm}p{2.0cm}p{1.0cm}p{2.0cm}p{1.0cm}p{2.0cm}p{1.0cm}}\hline
		$h$    &$\|\nabla e_{Q}^n\|_{L^{\infty}(L^2)}$    &Rate
		&$\|e_{Q}^n\|_{L^{\infty}(L^2)}$    &Rate
		&$|e_s^n|_{L^{\infty}(0,T)}$   &Rate   \\ \hline
		$1/32$    &1.46E-1     &---     & 6.10E-3  & ---    &4.71E-5        &---\\
		$1/64$    &3.59E-2     &2.03      &1.40E-3     &2.12    &1.34E-5         &1.81  \\
		$1/128$    &9.00E-3     &2.00    &3.49E-4     &2.02     &3.29E-6        &2.03  \\
		$1/256$    &2.30E-3    &2.00     & 8.70E-5     &2.00    &8.19E-7         &2.01   \\
		$1/512$    &5.64E-4   &2.00    &  2.17E-5     &2.00    &2.05E-7        &2.00   \\
		\hline
	\end{tabular}}
\end{table}



\begin{table}[!t]
	\renewcommand{\arraystretch}{1.1}
	\small
	\centering
	\caption{Error and convergence rates \textcolor{black}{of temporal discretization} for \textcolor{black}{\textnormal{MBP-sESAV2}} scheme \eqref{Q43}.
	}\label{table4}
	\scalebox{0.86}{\begin{tabular}{p{1.0cm}p{2.0cm}p{1.0cm}p{2.0cm}p{1.0cm}p{2.0cm}p{1.0cm}}\hline
		$\tau$    &$\|\nabla e_{Q}^n\|_{L^{\infty}(L^2)}$    &Rate
		&$\|e_{Q}^n\|_{L^{\infty}(L^2)}$    &Rate
		&$|e_s^n|_{L^{\infty}(0,T)}$   &Rate   \\ \hline
		$1/32$    &4.52E-4     &---     & 2.47E-5  & ---    &3.71E-5        &---\\
		$1/64$    &1.18E-4     &1.94     &6.43E-6     &1.94    &9.23E-6         &2.01  \\
		$1/128$    &3.00E-5     &1.97    &1.64E-6     &1.97     &2.29E-6        &2.01  \\
		$1/256$    &7.58E-6    &1.99     & 4.14E-7     &1.99    &5.72E-7         &2.00   \\
		$1/512$    &1.90E-6   &1.99    &  1.04E-7     &1.99    &1.43E-7        &2.00  \\
		\hline
	\end{tabular}}
\end{table}

\subsection{Disappearing hole}\label{holdd}
We use the value of the parameters and the initial data, as follows:
\bry
	&L_d=2,\quad\quad	L = \textcolor{black}{4.5}\times 10^{-3},\quad\quad
	a = -\textcolor{black}{4},\quad
	\quad c = \textcolor{black}{4},\\
		&Q_0=\frac{\textbf{n}_0\textbf{n}_0^T}{\|\textbf{n}_0\|_{h}^2}-\frac{1}{2}\textbf{I},\quad
		\textbf{n}_0=(\textcolor{black}{\frac{1}{16}}x(\textcolor{black}{2}-x)y(\textcolor{black}{2}-y),\
		\sin(\pi x)\sin(\pi y))^T.
\ery
\begin{figure}[h!]
	\centerline{\includegraphics[scale=0.075]{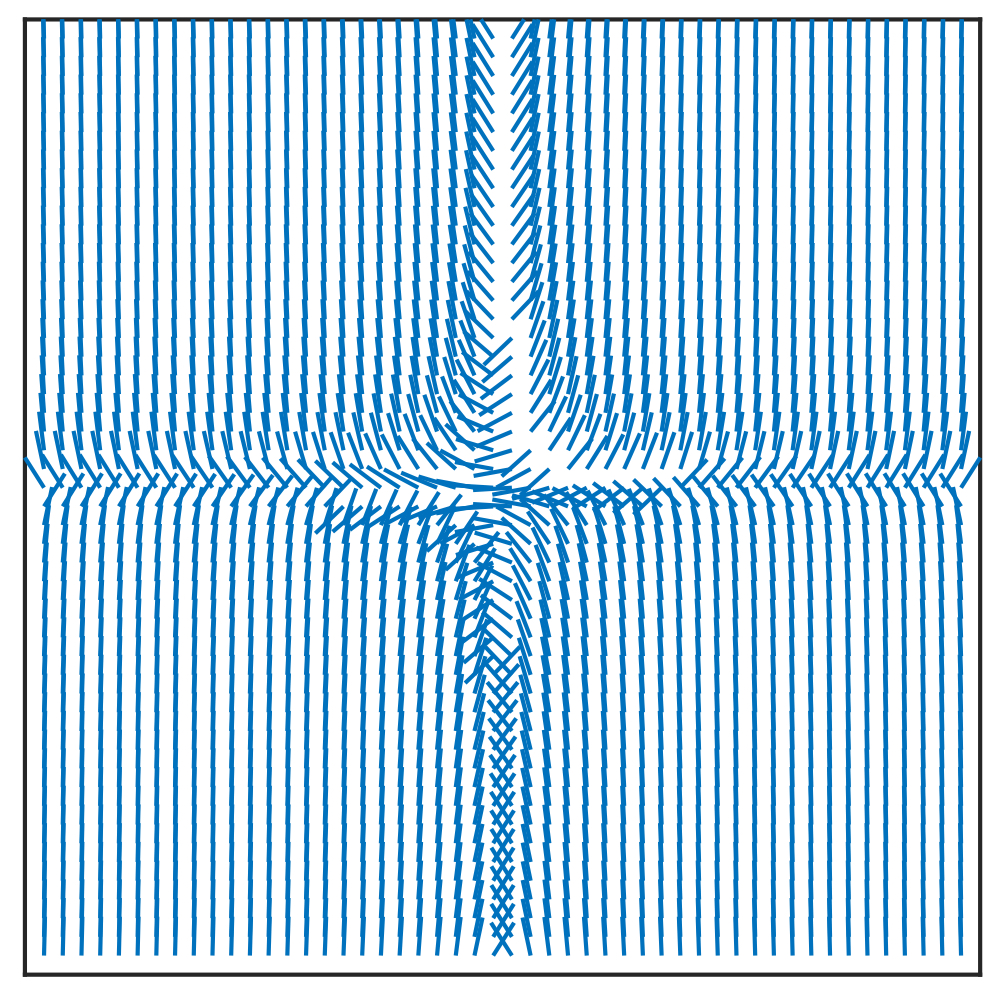}\includegraphics[scale=0.075]{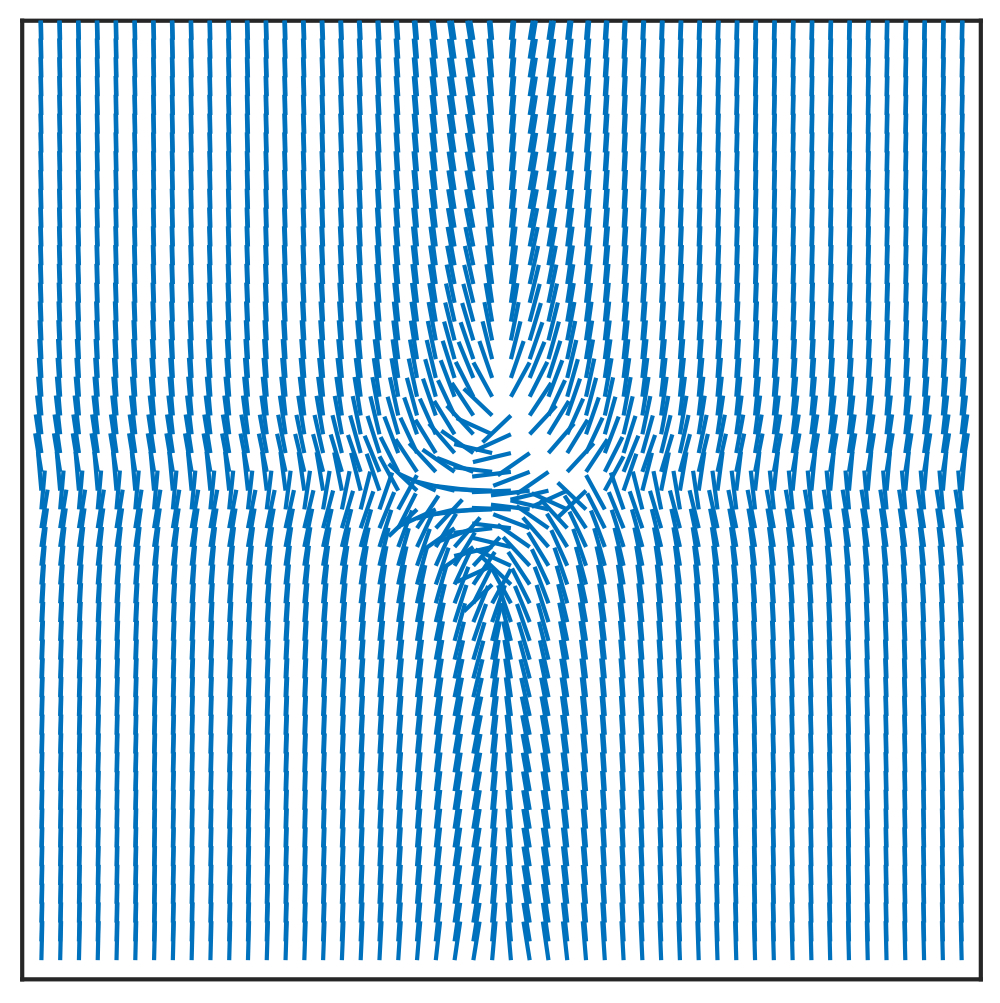}\includegraphics[scale=0.075]{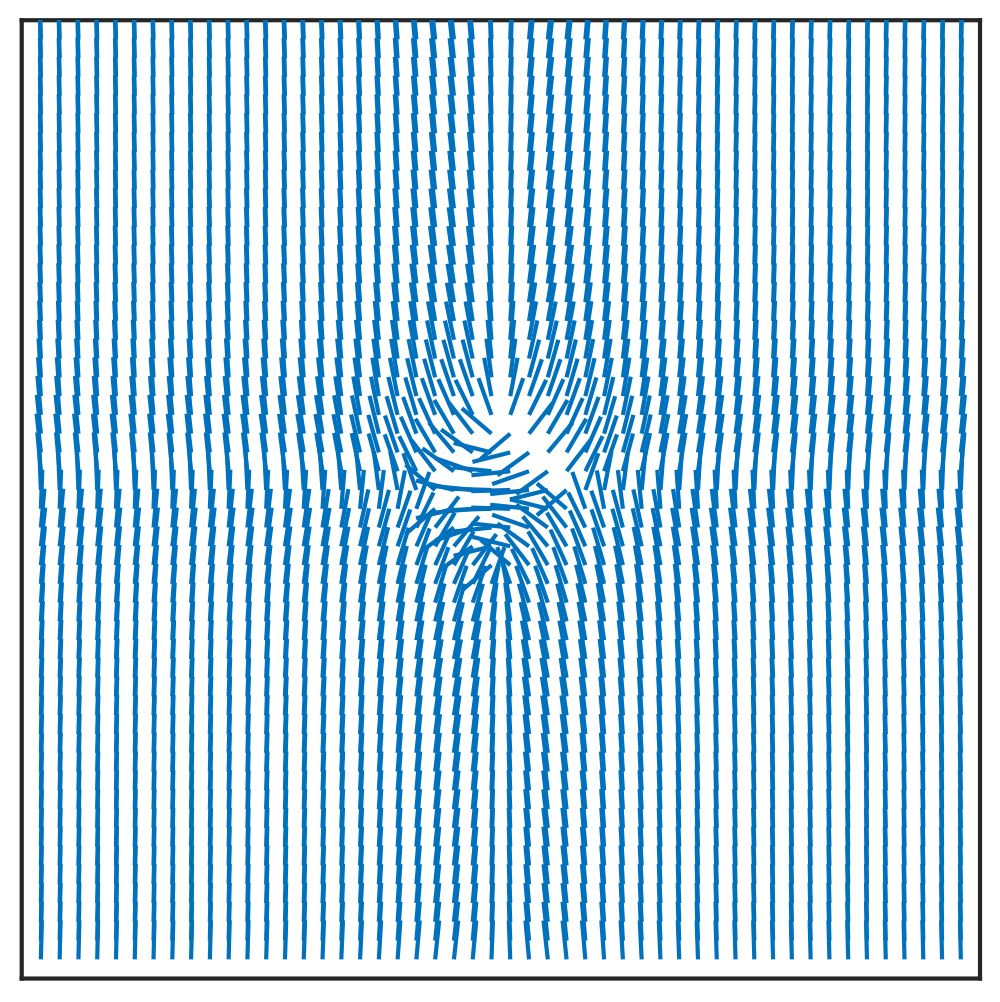}}
	\centerline{\includegraphics[scale=0.075]{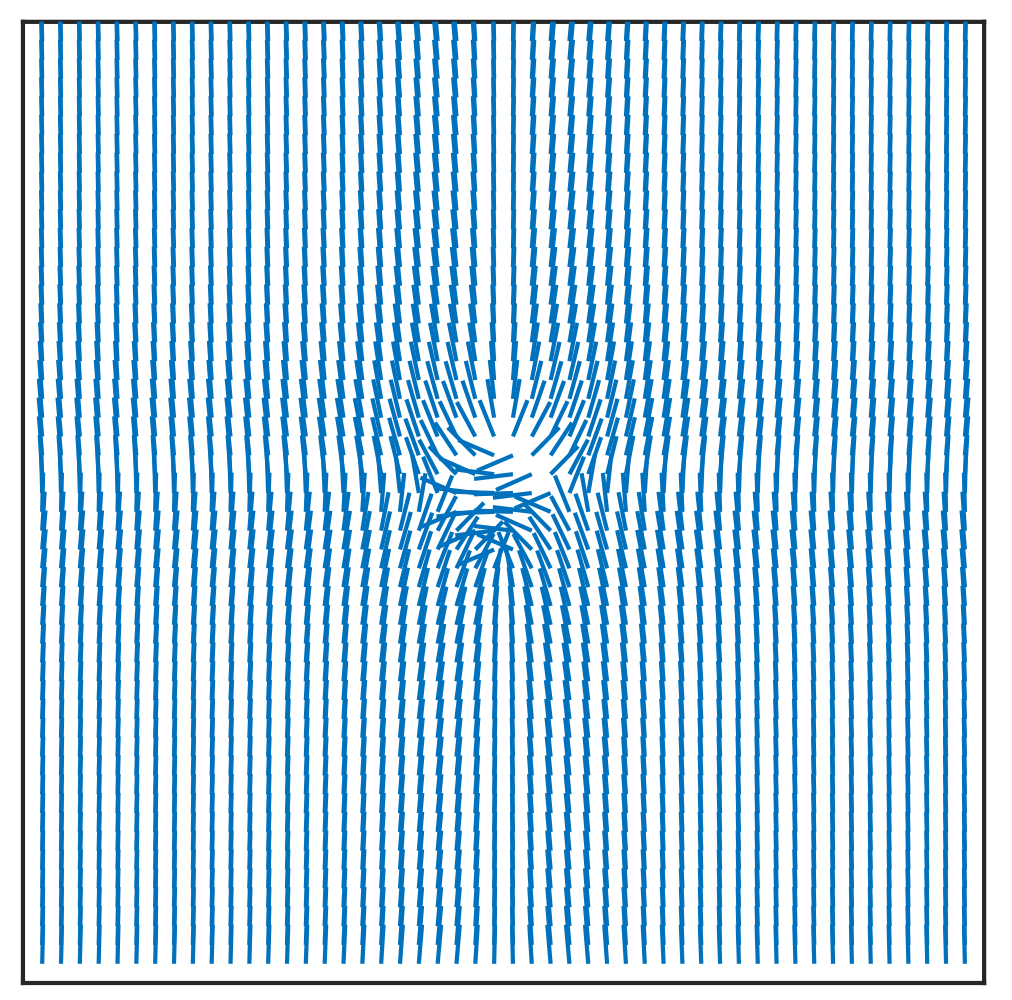}\includegraphics[scale=0.075]{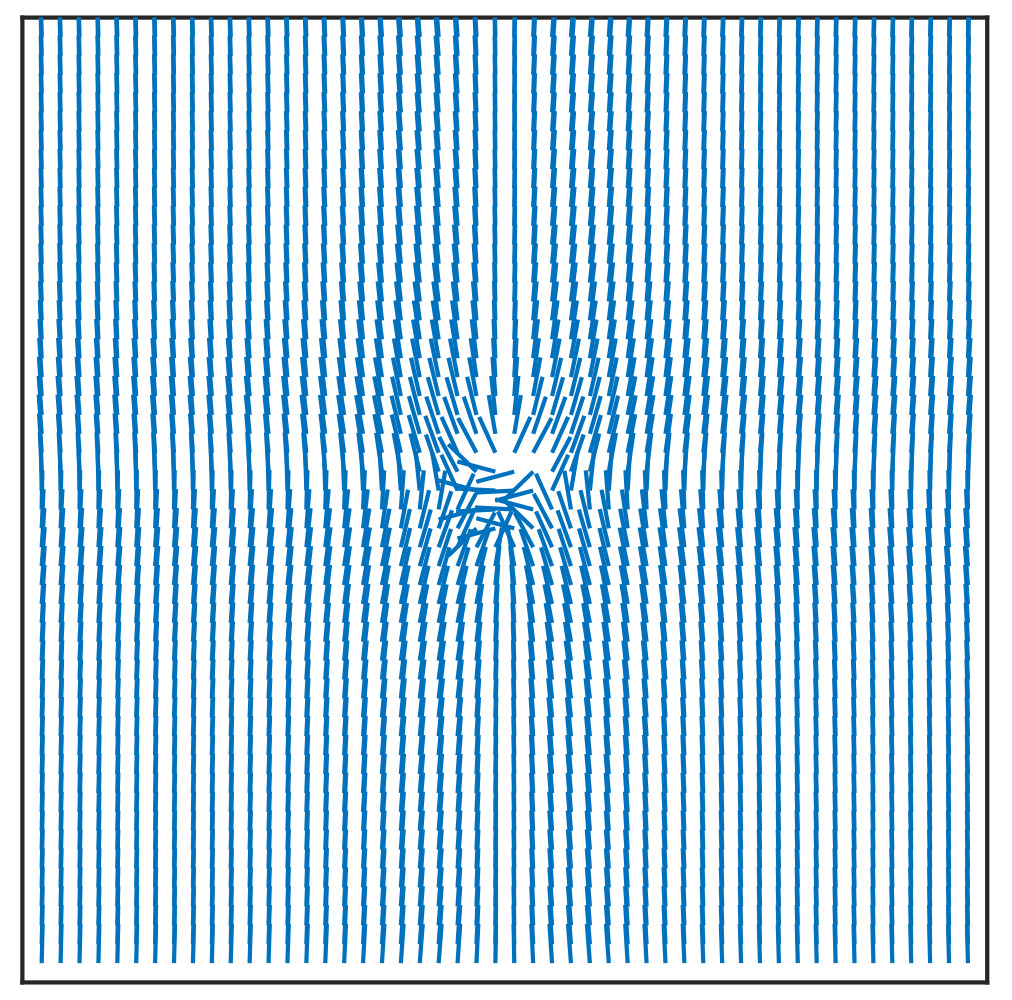}\includegraphics[scale=0.075]{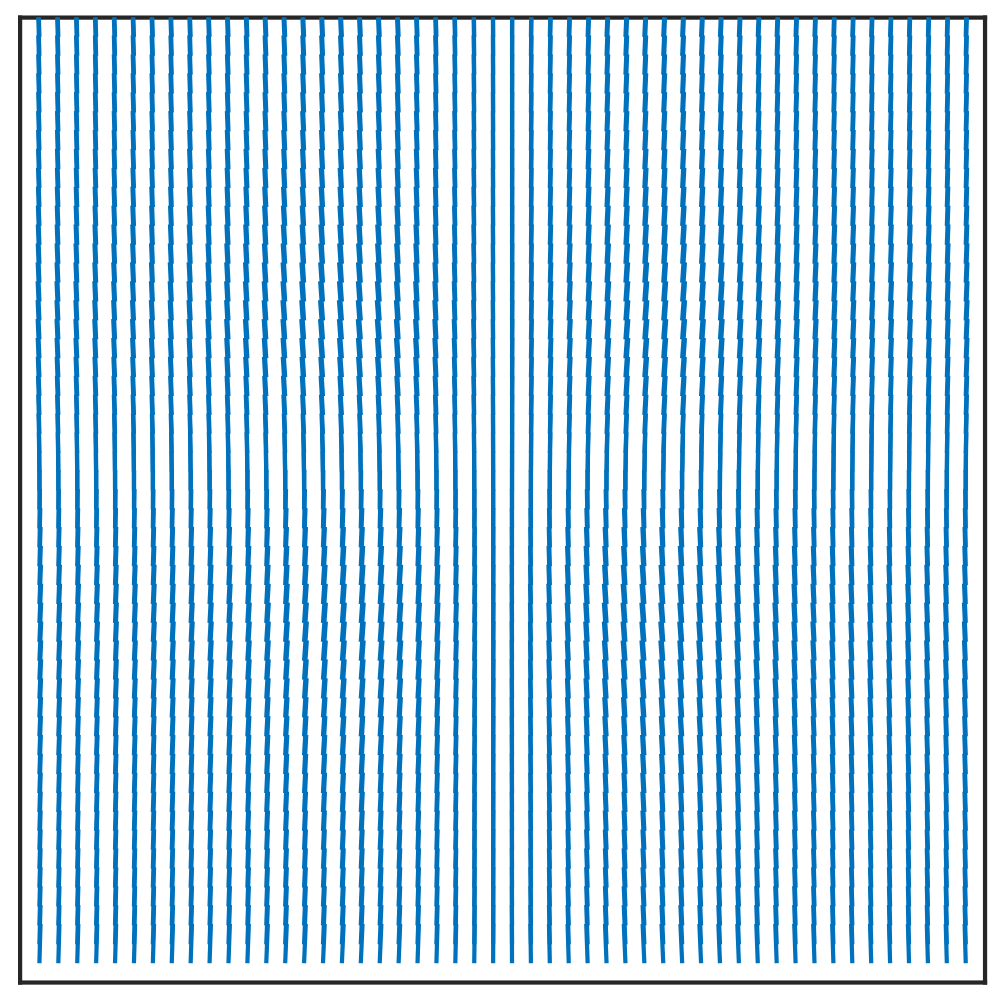}}
	\centerline{\includegraphics[scale=0.075]{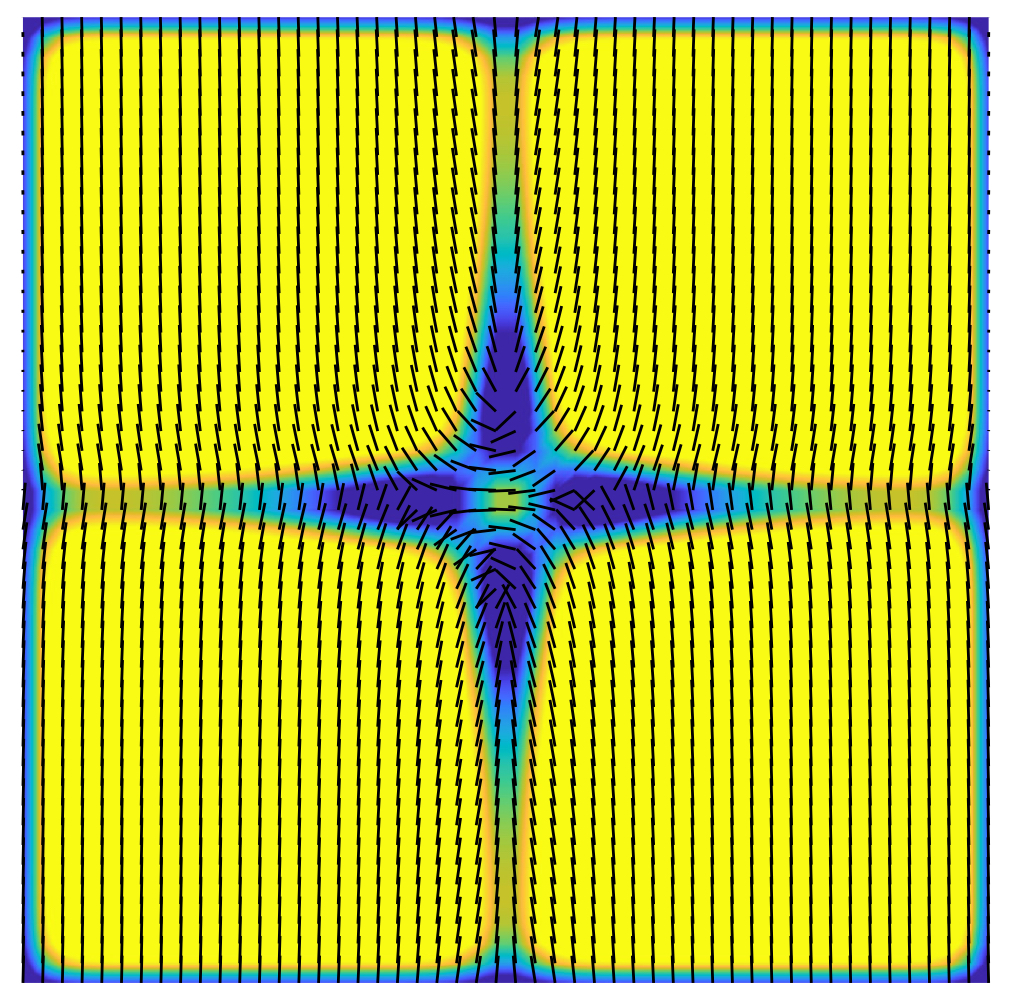}\includegraphics[scale=0.075]{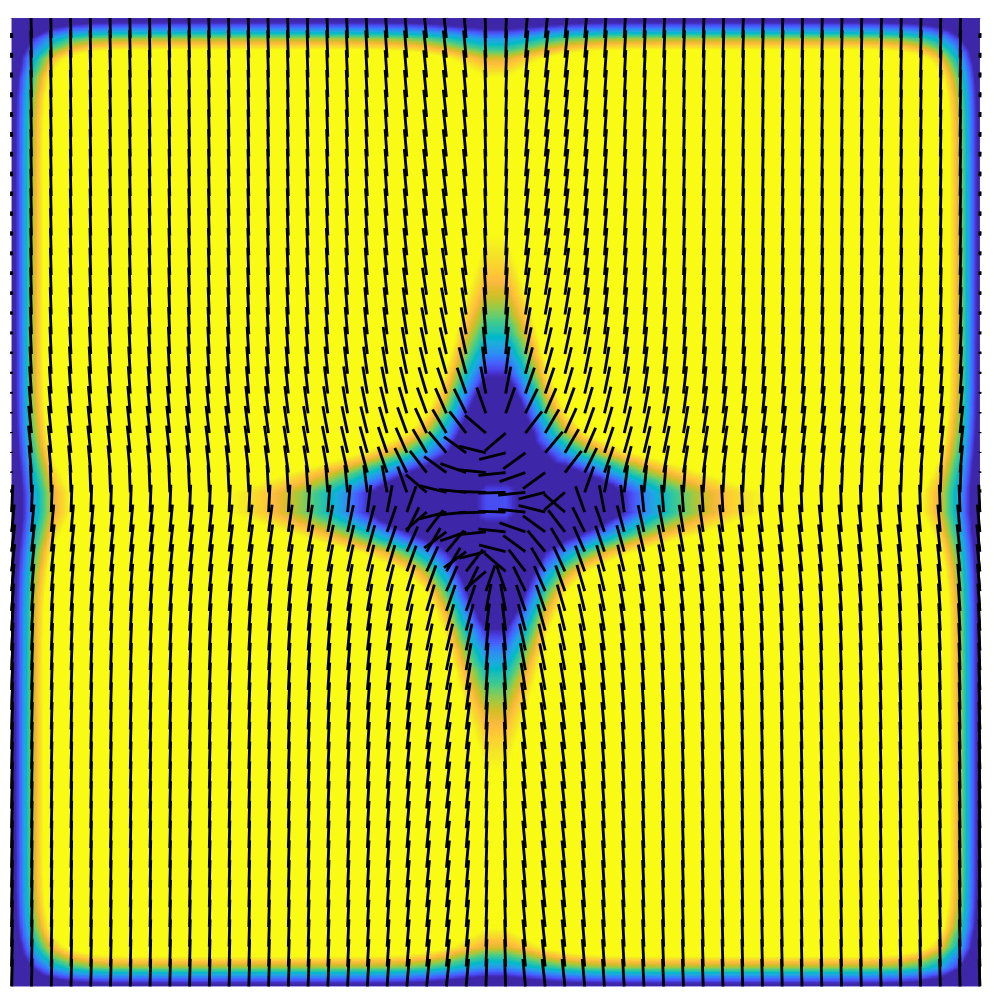}\includegraphics[scale=0.075]{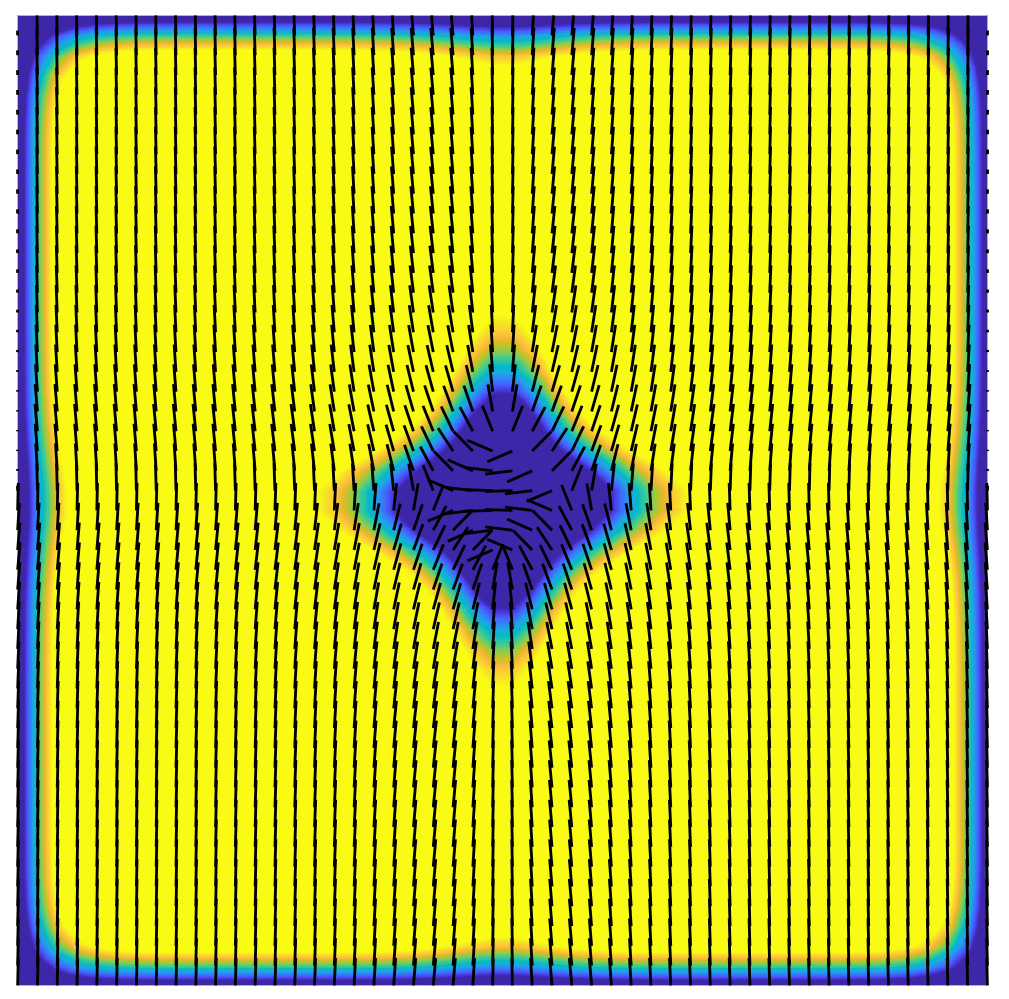}}
	\caption{Evolutions of disappearing holes at different time $t$.
		Snapshots are the major director orientation of the liquid crystal in the $xy$ plane taken at $t=0,\ 0.1,\ 0.2,\ 0.4,\ 0.7, \ 2.0$, respectively (top two rows).
		The difference of eigenvalues on the $xy$ plane for $Q+ \frac{1}{2}I$ at time $t = 0.1,\ 0.2,\ 0.4$, respectively (bottom row).} \label{hole_fig}
\end{figure}

\begin{figure*}[h!]
\begin{minipage}[t]{0.325\linewidth}
\centerline{\includegraphics[scale=0.245]{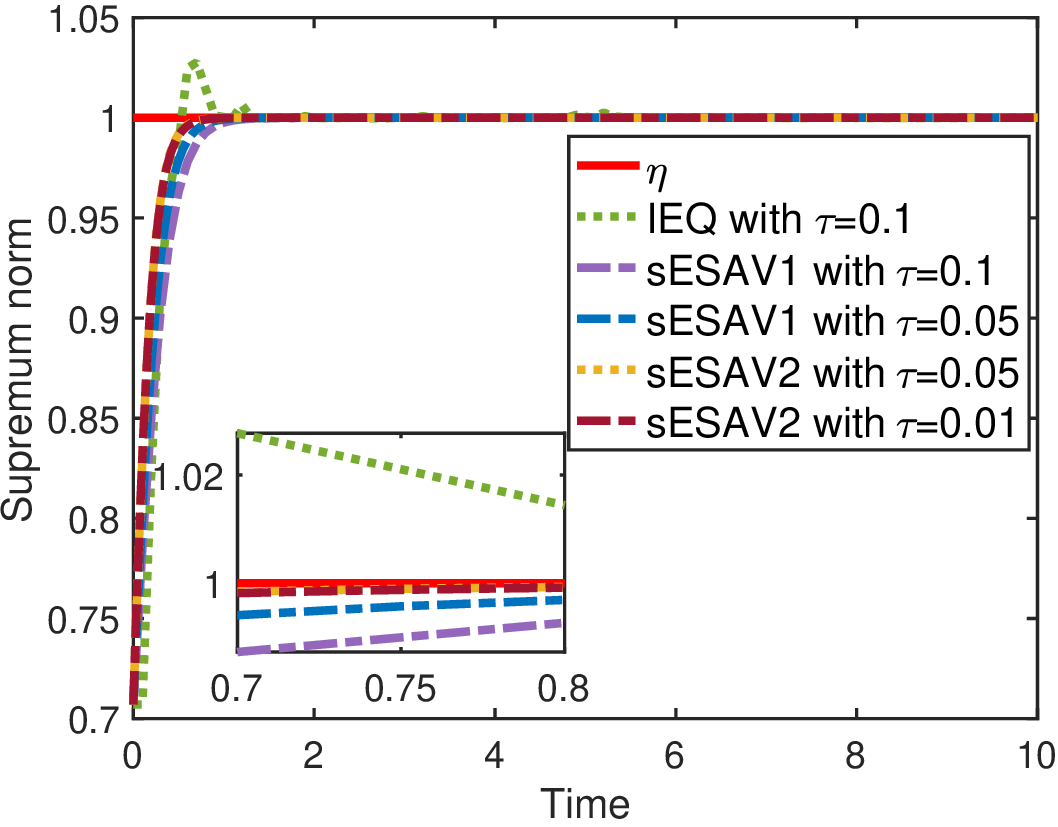}}
\centerline{\footnotesize{(a)} }
\end{minipage}
\begin{minipage}[t]{0.325\linewidth}
\centerline{\includegraphics[scale=0.245]{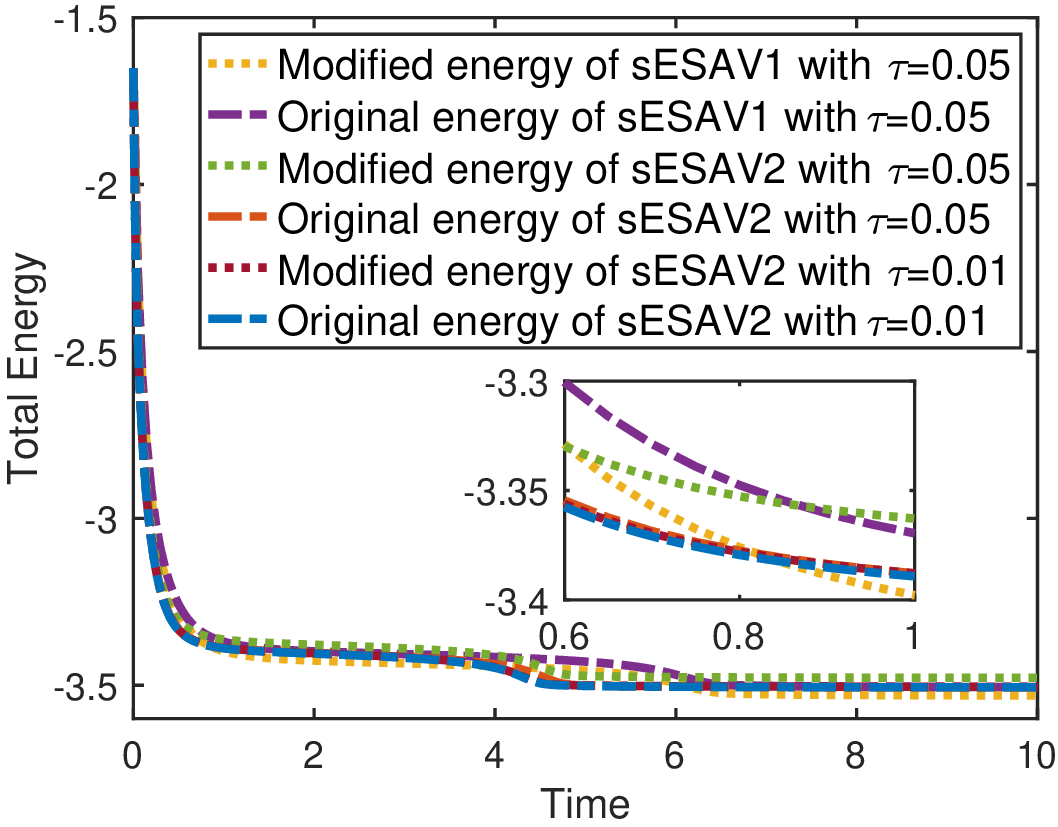}}
\centerline{\footnotesize{(b)} }
\end{minipage}
\begin{minipage}[t]{0.325\linewidth}
\centerline{\includegraphics[scale=0.245]{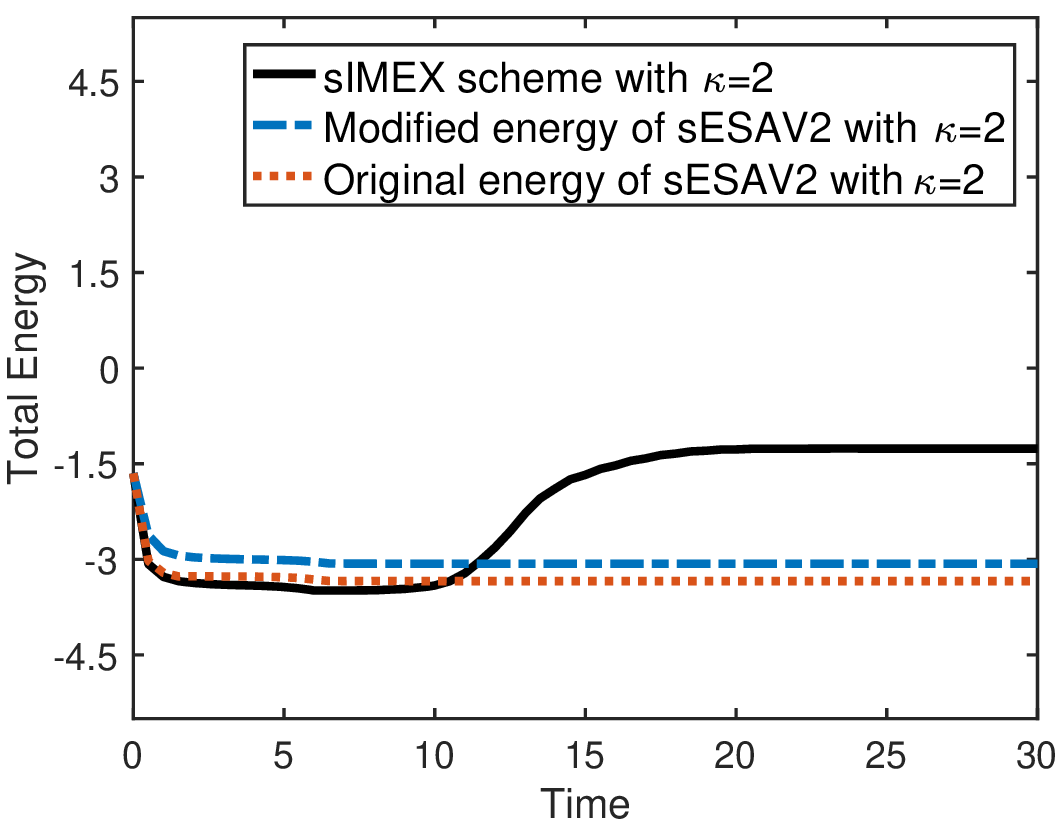}}
\centerline{\footnotesize{(c)}}
\end{minipage}
\vspace{-0.3cm}
\caption{Evolutions of supremum norms \textcolor{black}{ \textnormal{(a)}} and energies \textcolor{black}{ \textnormal{(b)}} of simulated solutions. \textcolor{black}{ \textnormal{(c)} Comparison of the sIMEX scheme with \textnormal{sESAV2} scheme using $\tau=0.5$ and $\kappa=2$.
		}}\label{hole}
\end{figure*}

In the computations, we choose
\textcolor{black}{$\kappa=8$ and  $80$ grid points in space in each dimension}.
Figure \ref{hole_fig} presents the major director  orientation of the liquid crystal in the $xy$ plane.
It can be observed that the initial misalignment vanishes
along the axes and then propagates in a shrinking circle toward the center of the domain and disappears eventually.
Figure \ref{hole_fig} shows that the difference between the two eigenvalues on the $xy$ plane is a true measure of the orientation of liquid crystals. When the difference value is zero, it is called a defect. It represents a state of liquid crystal, indicating that the orientation on the $xy$ plane is isotropic and dominant.
The comparison of the supremum norms of the simulated solutions based on the IEQ scheme and the MBP-sESAV approach \textcolor{black}{with different time steps} are presented in \textcolor{black}{Figure \ref{hole}(a)}.
It shows that our MBP-sESAV  scheme performs better than the IEQ method presented in \cite{gudibanda2022convergence} with the same time step.
The evolutions of the energies of simulated solutions are shown in Figure \ref{hole}(b) with \textcolor{black}{MBP-sESAV1 \eqref{Q24} and MBP-sESAV2 \eqref{Q46}}.
Figures \ref{hole}(a)-\ref{hole}(b) shows \textcolor{black}{MBP-sESAV1} preserves the MBP and the energy dissipation perfectly even with \textcolor{black}{larger} time steps.
\textcolor{black}{
Moreover, we observe that the stabilized implicit-explicit (sIMEX) scheme leads to an increase in energy, while the sESAV2 scheme exhibits energy dissipation in Figure \ref{hole}(c). 
}

\subsection{Dynamics of orientation in the liquid crystal}
In this subsection, \textcolor{black}{an} example \textcolor{black}{is} performed to show orientation dynamics  in flows of liquid crystals numerically using the implemented linear scheme.
\textcolor{black}{The space is discretized with 100 grid points in each dimension.}
Here, we adopt the variable time step sizes $\tau_{n+1}$ updated
by using the approach from \cite{qiao2011adaptive}
\begin{equation*}
	\tau_{n+1}=\max
	\{ \tau_{\min},\frac{\tau_{\max}}{\sqrt{1+\alpha|d_t\mathcal{E}_h[Q_{h}^n]|^2}}\},
\end{equation*}
where $d_t\mathcal{E}_h[Q_{h}^n]=(\mathcal{E}_h[Q_{h}^{n+1}]-\mathcal{E}_h[Q_{h}^n])/\tau_n$
and $\alpha > 0$ is a constant parameter.

In this example, we simulate the orientation dynamics in 3D liquid crystal flow.
For the parameter values, we set
\begin{equation*}
	L_d=2,\quad	L = 1.0\times 10^{-3},\quad
	a = -1.25,
	\quad b = 0.25,
	\quad c = 1,
	\quad \kappa=6.
\end{equation*}
The initial state for $Q_0$ is defined as $Q_0=\textbf{n}_0\textbf{n}_0^T/\|\textbf{n}_0\|_{h}^2-\frac{1}{3}\textbf{I}$
with
\begin{equation*}
	\textbf{n}_0=\left\{
	\begin{array}{rcl}
		(1,0,0)^T, &  & x\in[1.15,1.65]\times[0.75,1.25]\times[0.35,0.85]; \\
		(0,0,1)^T, &  & x\in[0.35,0.85]\times[0.75,1.25]\times[1.15,1.65]; \\
		(0,1,0)^T, &  & \text{otherwise}.
	\end{array}
\right.
\end{equation*}

\begin{figure}[t!]
	\centerline{\includegraphics[scale=0.08]{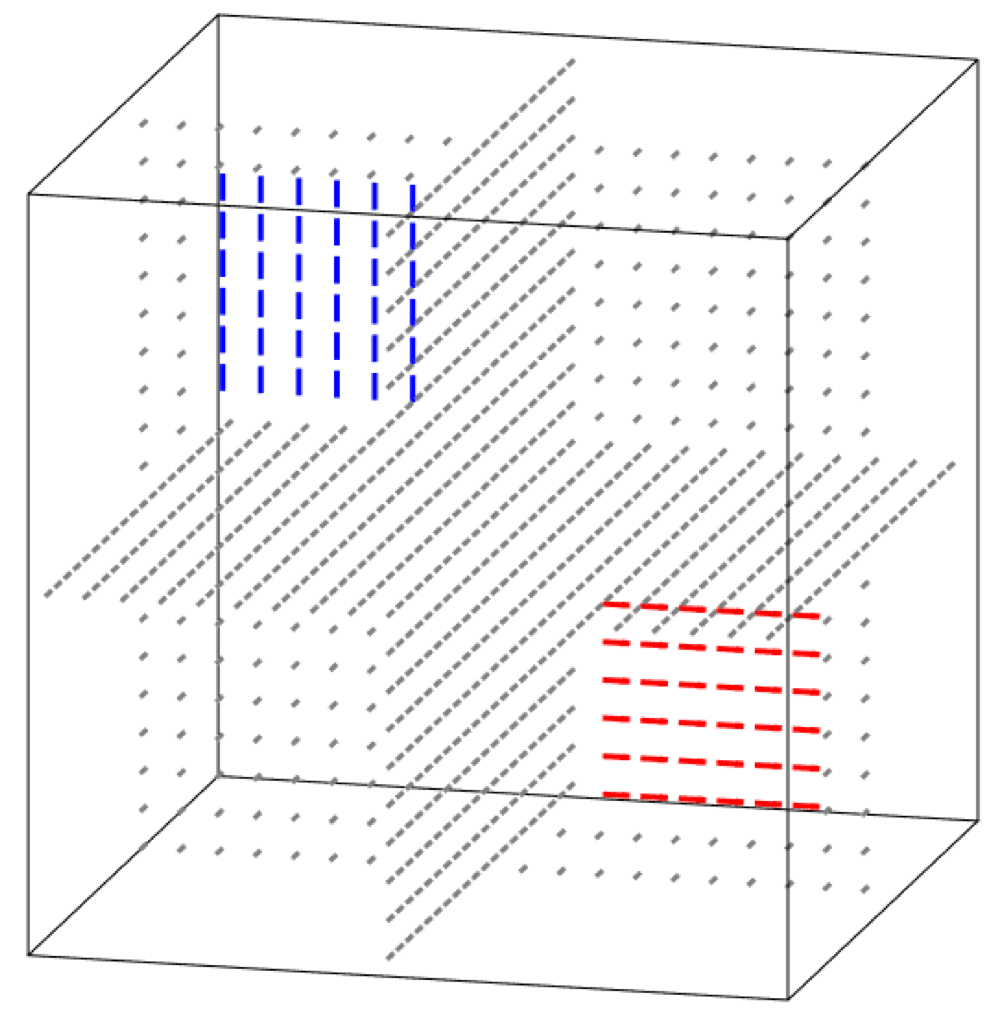}\includegraphics[scale=0.08]{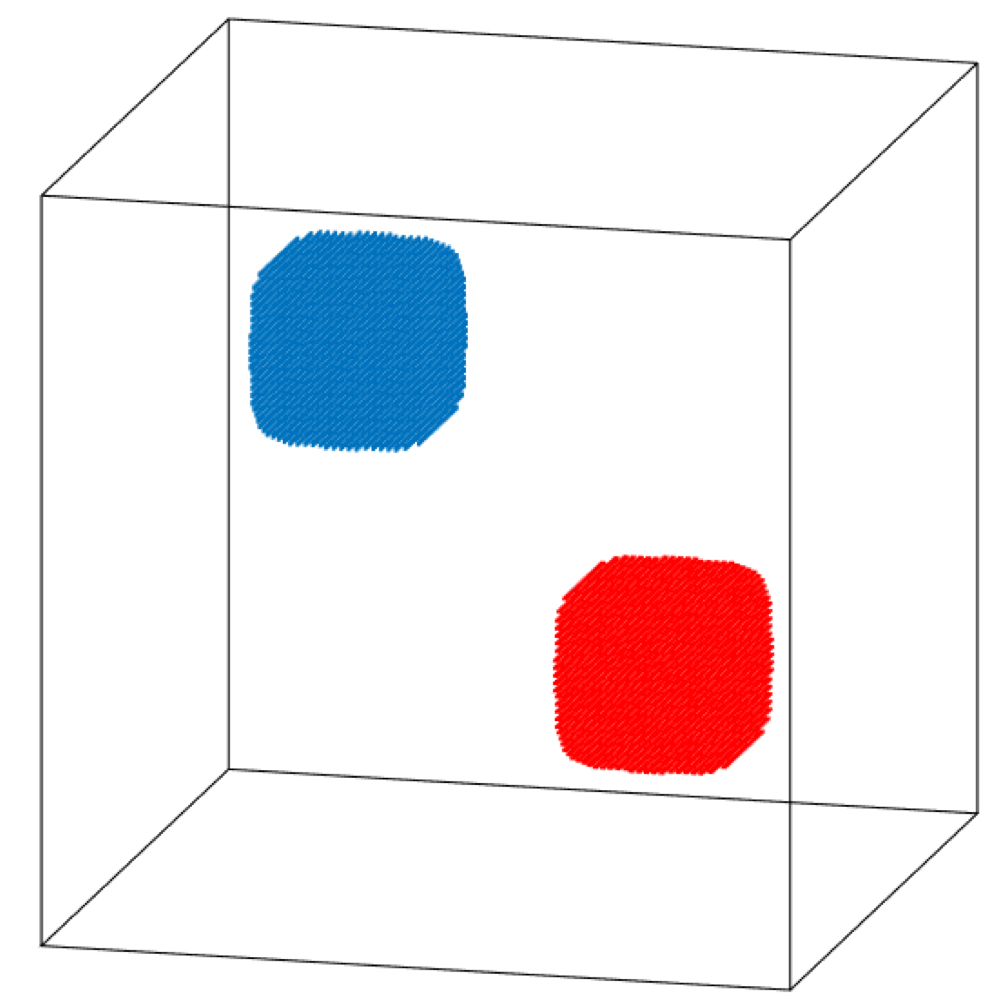}\includegraphics[scale=0.08]{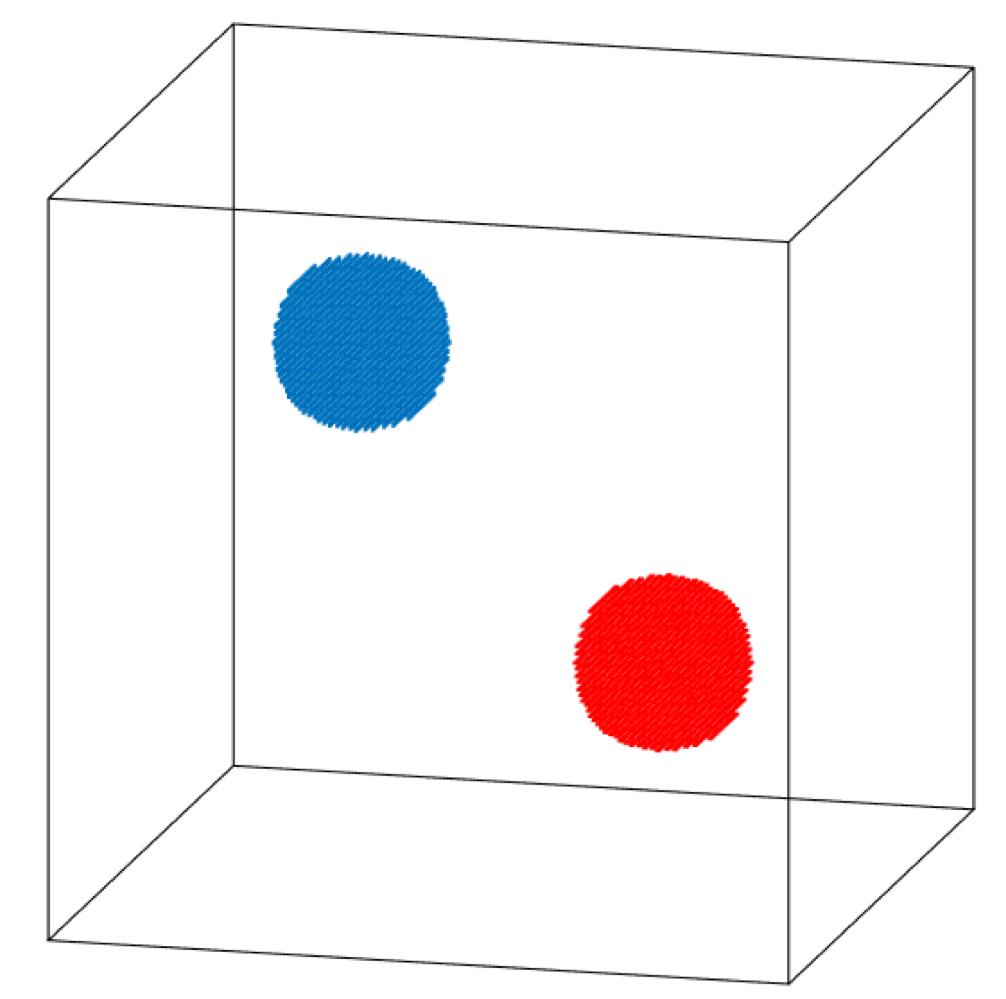}}
	\centerline{\includegraphics[scale=0.08]{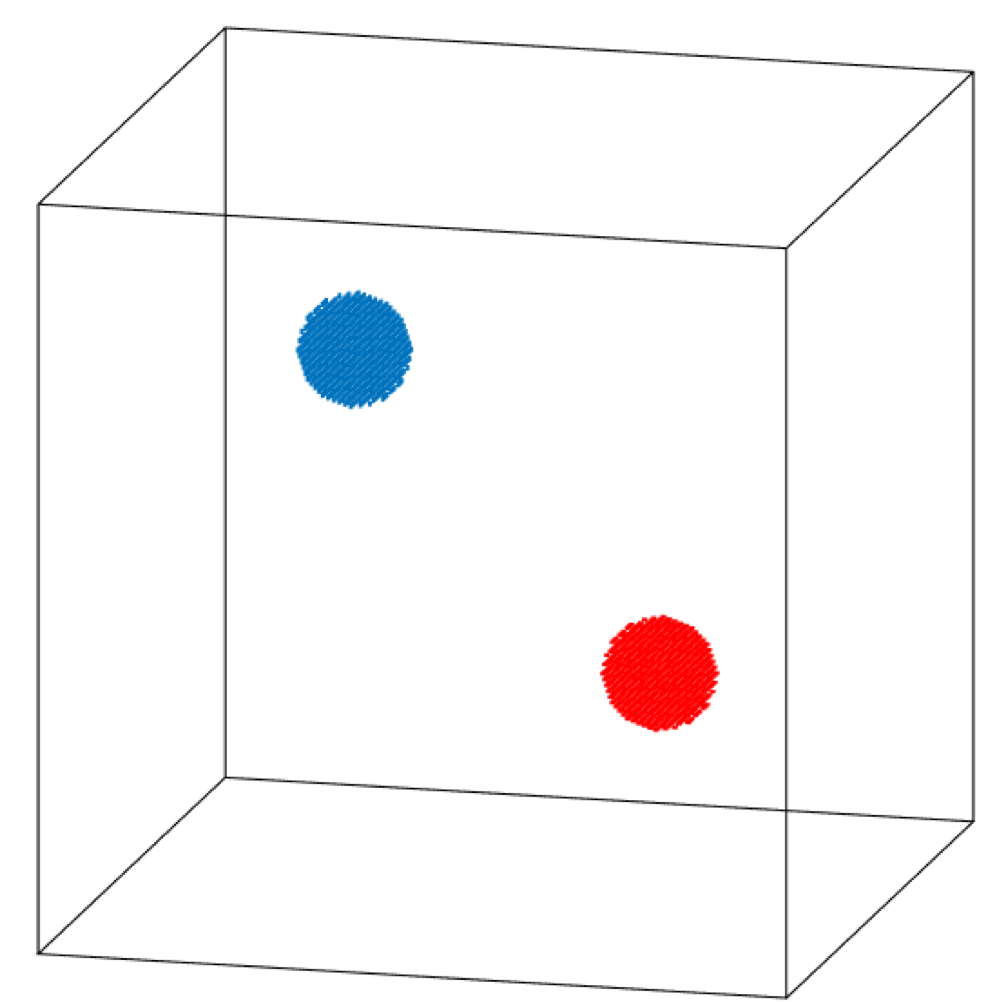}\includegraphics[scale=0.08]{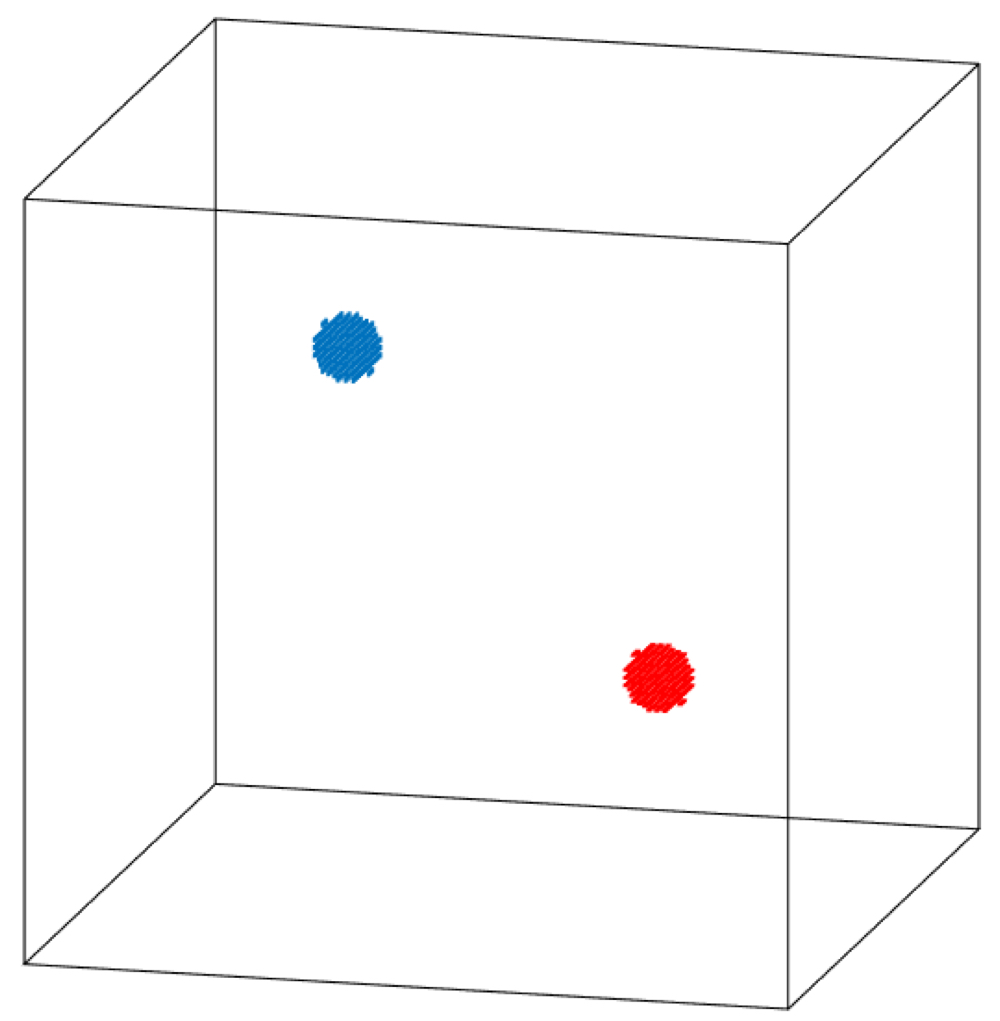}\includegraphics[scale=0.08]{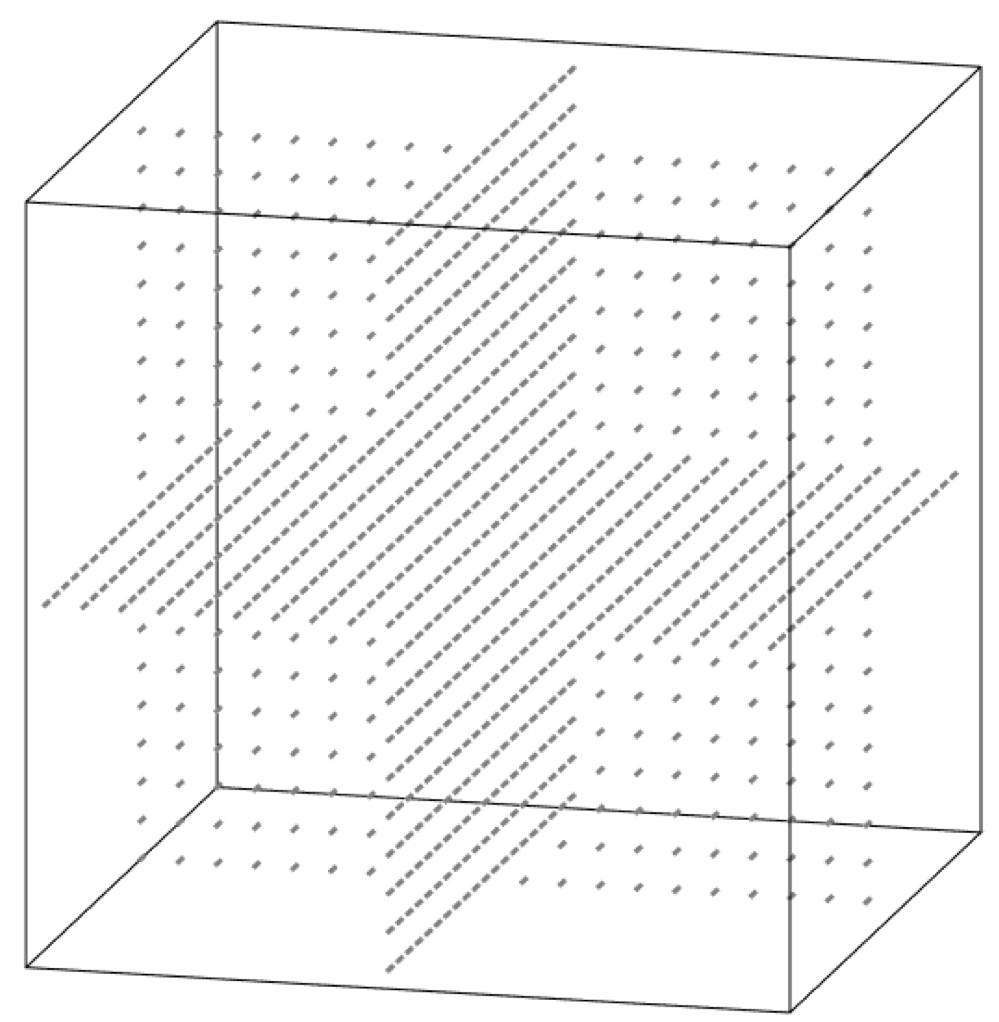}}
	\centerline{\includegraphics[scale=0.08]{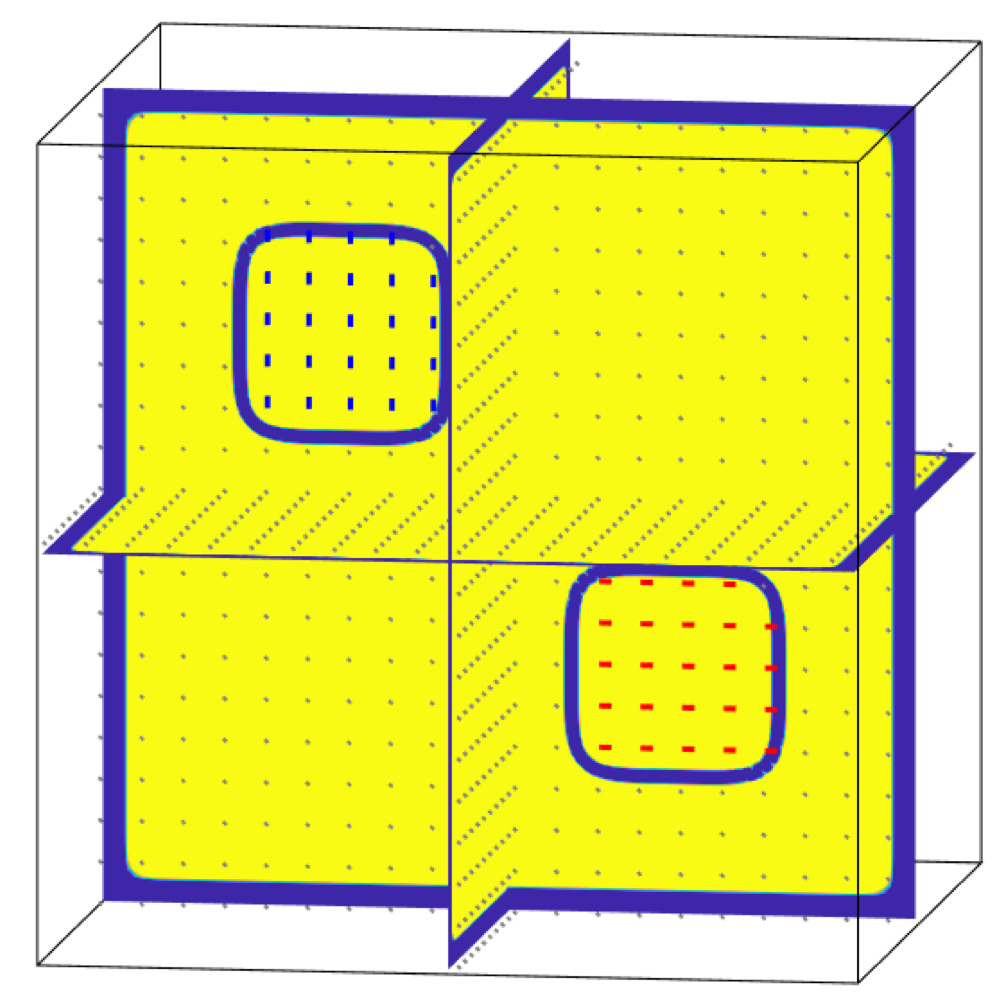}\includegraphics[scale=0.08]{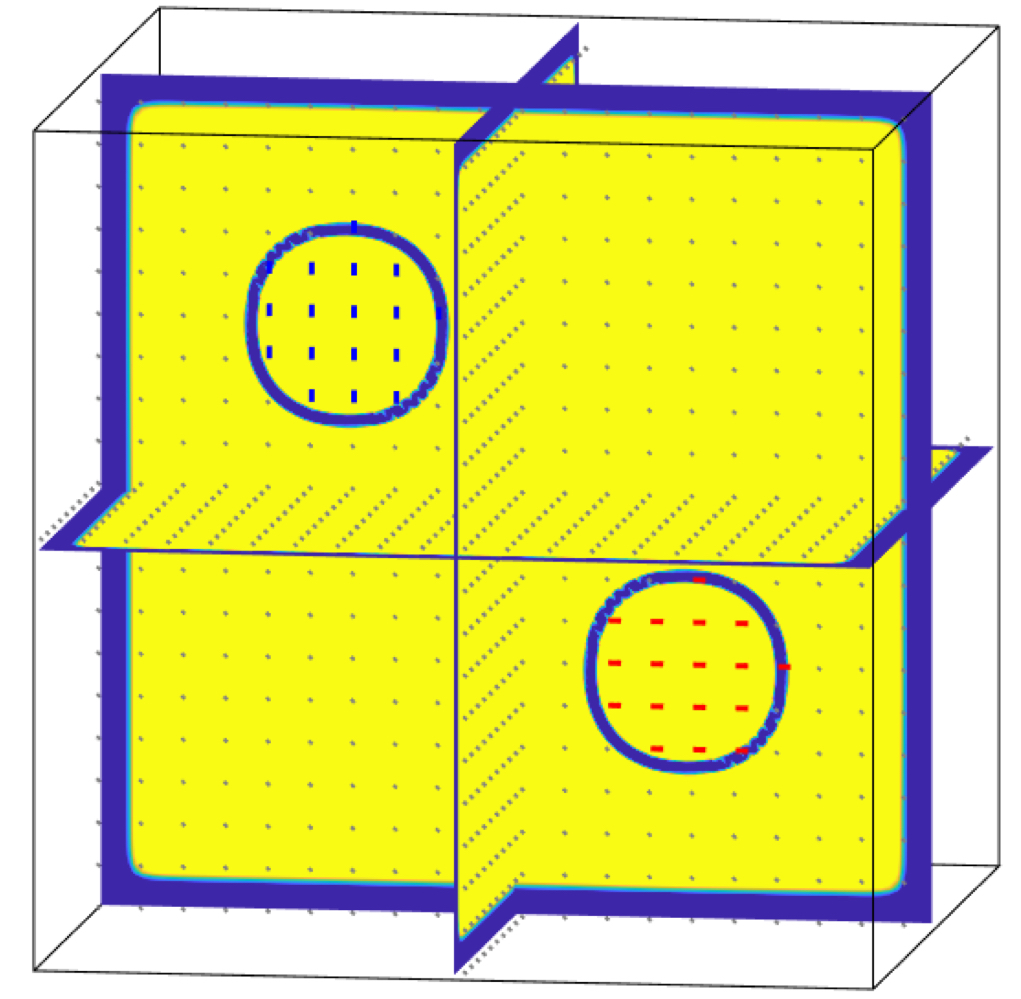}\includegraphics[scale=0.08]{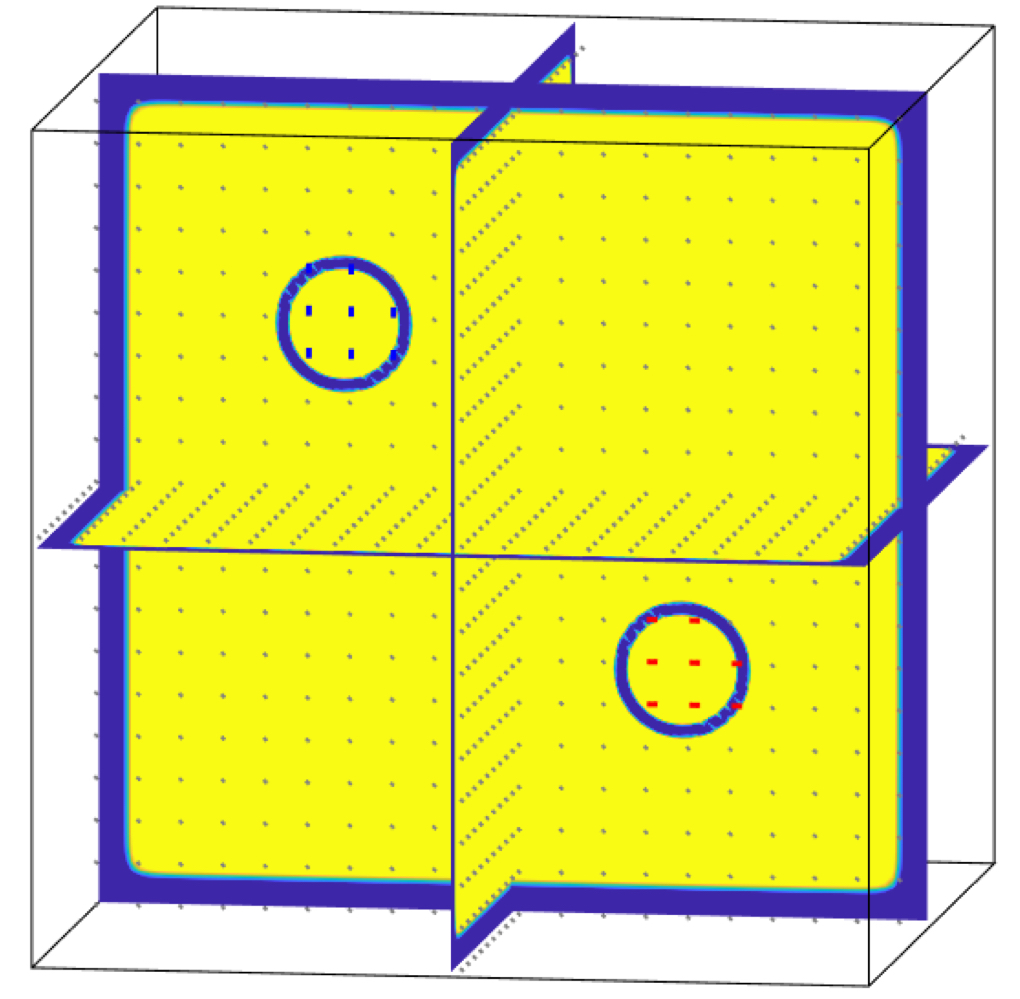}}
	\caption{Orientation of liquid crystal at different time $t$.
		Snapshots are the major director orientation of the liquid crystal in the $xy$ plane taken at $t=0,\ 2,\ 10,\ 18,\ 21, \ 30$, respectively (top two rows).
		The difference of eigenvalues on the $xyz$ plane for $Q+ \frac{1}{3}I$ at time $t = 2,\ 10,\ 18$, respectively (bottom row).}\label{circle3D}
\end{figure}

\begin{figure}[t!]
	\centerline{\includegraphics[scale=0.24]{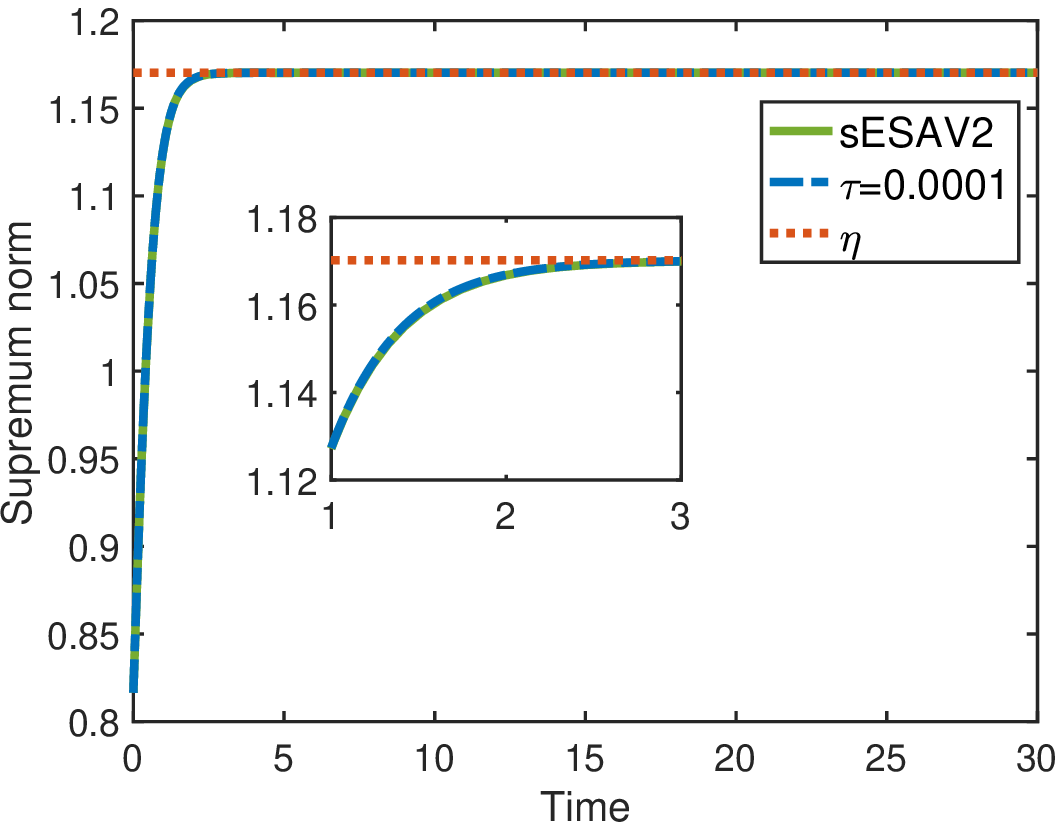}
	\includegraphics[scale=0.24]{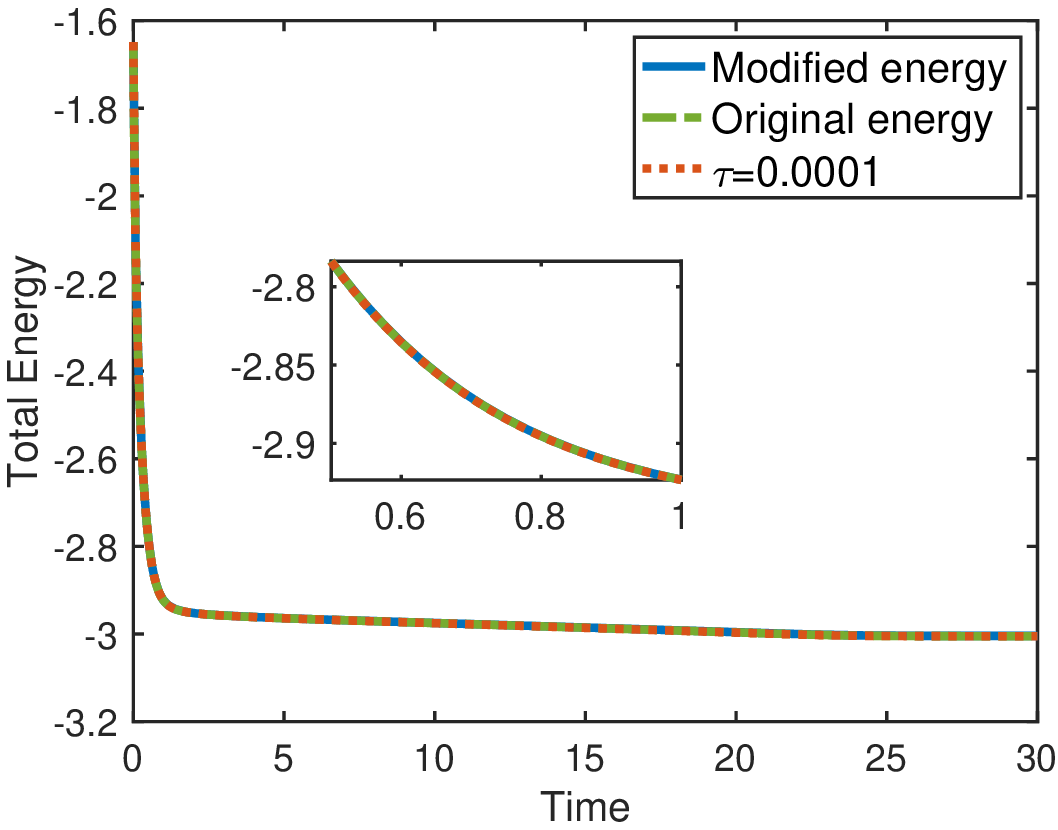}
		\includegraphics[scale=0.24]{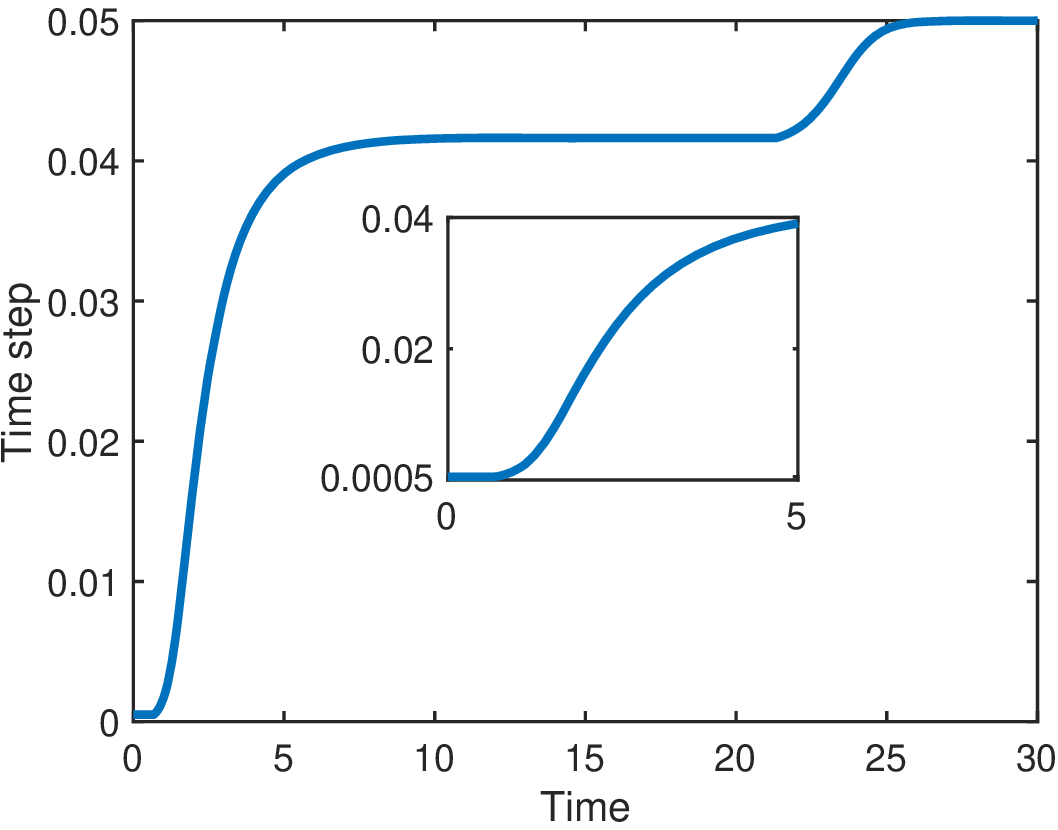}}
	\caption{Evolutions of the supremum norms (left) and the energies \textcolor{black}{(middle) and the adaptive time steps (right).} }\label{circle3}
\end{figure}
\textcolor{black}{
In this example, we set the minimal and maximal time step sizes as $\tau_{min} = 0.0005$ and $\tau_{max} = 0.05$,
respectively, with $\alpha = 10^5$ being defined.}
The orientation of liquid crystal on the $xy$ plane, which is denoted by $\textbf{n}$, is mainly described by the dominant eigenvalue  of $Q$.
For the 3D case, the orientation of the liquid crystal is depicted in Figures \ref{circle3D}. 
Apparently, there are more gray molecules than black and red molecules at $T = 0$ in Figure \ref{circle3D}.
The set with red and black molecules gets smaller with time,
then the defects stay in a spherical form and the spherical shrinks with time.
The liquid crystal directions get closer to achieving a uniform vertical shape finally, i.e., the spherical disappears.

The evolutions of the supremum norm of $|Q|_F$ and the energy of simulated solutions are displayed in \textcolor{black}{Figure \ref{circle3}}.
\textcolor{black}{We obtained the solution by implementing the MBP-sESAV2 scheme with adaptive time stepping and an exceptionally small time step of  $\tau= 0.0001$.}
	Similar to \textcolor{black}{Section \ref{holdd}}, the sESAV2 schemes preserve the MBP and energy dissipation law as expected.
Moreover, it can also be verified that the maximum bound of the numerical scheme we obtained is bounded by the theoretical bound $\eta$.

\section{Concluding remarks}
In this paper, we proposed two linear and efficient schemes for the $Q$-tensor flow of liquid crystals based on the sESAV approach.
The constructed schemes are unconditionally energy-stable and preserve the MBP. Moreover, we establish
rigorous error estimates for the second-order scheme. One of the future research directions is to extend the approaches in this paper to hydrodynamic $Q$-tensor models, which couples a Navier-Stokes system for the fluid velocity with a parabolic reaction-convection-diffusion system.

\bibliographystyle{siamplain}
\bibliography{tes}

\end{document}